\documentclass[11pt,a4paper]{article}%
\usepackage{amsfonts}
\usepackage{geometry}
\usepackage{amsmath}
\usepackage{amssymb}
\usepackage{amsthm}
\usepackage{hyperref}
\usepackage{mathtools}
\usepackage{graphicx}
\usepackage[utf8]{inputenc}
\usepackage[T1]{fontenc}
\usepackage[document]{ragged2e}
\usepackage{placeins}
\usepackage{enumitem}
\usepackage{float}
\usepackage{comment}
\usepackage{hhline}
\usepackage{tikz}
\usepackage{parskip}
\providecommand{\U}[1]{\protect\rule{.1in}{.1in}}
\usetikzlibrary{calc}
\usetikzlibrary{patterns}
\usepackage[
backend=bibtex,
autocite=inline,
autolang=hyphen,
autopunct,
style=numeric,
sorting=nty,
giveninits=true,
doi=false,
url=false,
isbn=false
]{biblatex}
\AtEveryBibitem{\clearfield{pages}} 
\AtEveryBibitem{\clearfield{note}} 
\numberwithin{equation}{section}

\newtheorem{theorem}{Theorem}[section]
\newtheorem{corollary}[theorem]{Corollary}
\newtheorem{definition}[theorem]{Definition}
\newtheorem{lemma}[theorem]{Lemma}
\newtheorem{remark}[theorem]{Remark}
\newtheorem{example}[theorem]{Example}
\newtheorem{proposition}[theorem]{Proposition}

\newtheorem{assumption}[theorem]{Assumption}

\makeatletter
\def\moverlay{\mathpalette\mov@rlay}
\def\mov@rlay#1#2{\leavevmode\vtop{   \baselineskip\z@skip \lineskiplimit-\maxdimen
\ialign{\hfil$\m@th#1##$\hfil\cr#2\crcr}}}
\newcommand{\charfusion}[3][\mathord]{
#1{\ifx#1\mathop\vphantom{#2}\fi
\mathpalette\mov@rlay{#2\cr#3}
}
\ifx#1\mathop\expandafter\displaylimits\fi}
\makeatother

\def\XXint#1#2#3{{\setbox0=\hbox{$#1{#2#3}{\int}$ }
\vcenter{\hbox{$#2#3$ }}\kern-.6\wd0}}

\newcommand{\vertiii}[1]{{\left\vert\kern-0.25ex\left\vert\kern-0.25ex\left\vert #1
\right\vert\kern-0.25ex\right\vert\kern-0.25ex\right\vert}}

\let\originalleft\left
\let\originalright\right
\renewcommand{\left}{\mathopen{}\mathclose\bgroup\originalleft}
\renewcommand{\right}{\aftergroup\egroup\originalright}

\begin{document}
\title{Pressure-improved Scott-Vogelius type elements}
\author{Nis-Erik Bohne\thanks{(nis-erik.bohne@math.uzh.ch), Institut f\"{u}r
Mathematik, Universit\"{a}t Z\"{u}rich, Winterthurerstr 190, CH-8057
Z\"{u}rich, Switzerland}
\and Benedikt Gr\"{a}\ss le\thanks{(graesslb@math.hu-berlin.de), Institut f\"{u}r
Mathematik, Humboldt-Universit\"{a}t zu Berlin, 10117 Berlin, Germany }
\and Stefan A. Sauter\thanks{(stas@math.uzh.ch), Institut f\"{u}r Mathematik,
Universit\"{a}t Z\"{u}rich, Winterthurerstr 190, CH-8057 Z\"{u}rich,
Switzerland}}
\maketitle
\justifying
\begin{abstract}
The Scott-Vogelius element is a popular finite element for the discretization
of the Stokes equations which enjoys inf-sup stability and gives
divergence-free velocity approximation. However, it is well known that the
convergence rates for the discrete pressure
deteriorate in the presence of certain \textit{critical vertices} in
a
triangulation of the domain. Modifications of the Scott-Vogelius element such
as the recently introduced pressure-wired Stokes element also suffer from this
effect. In this paper we introduce a simple modification strategy for these
pressure spaces that preserves the inf-sup stability while the pressure
converges at an optimal rate.
\end{abstract}




\section{Introduction}

In this paper we consider the numerical solution of the stationary Stokes
equations by conforming Galerkin finite element methods. This is a vivid
research topic since many decades in numerical analysis and scientific
computing. The unknowns are the vector-valued velocity field and the scalar
pressure and the Galerkin discretization is based on the choice of a pair of
finite element spaces: one, say $\mathbf{S}$, for the velocity and one, say
$M$, for the pressure approximation. It is well-known that the most
\textit{intuitive} choice for $\mathbf{S}$, i.e., continuous, piecewise
polynomials of degree $k$ and for $M$, i.e., discontinuous, piecewise
polynomials of degree $k-1$ (which we will denote as
the full pressure space) can be unstable: although the
continuous problem is well-posed the Galerkin discretization may result in a
singular system matrix and is not solvable (see, e.g.,
\cite{vogelius1983right} and~\cite[Chap.~7]{Braessengl} for
quadrilateral meshes).

This problem
motivated the development of many pairs $\left(  \mathbf{S},M\right)  $ of
Stokes elements; standard strategies include enrichment of the intuitive
velocity space, see, e.g., \cite{merdon_SV}, reducing the intuitive pressure
space, see, e.g., \cite{Bernardi_Raugel_1985}, \cite{GuzmanScott2019},
\cite{ScottVogelius}, \cite{vogelius1983right} and combinations thereof
\cite{Ainsworth_parker_hp_version}, \cite{arnold1984stable}, \cite[Chap. 3, \S 7]{Braessengl},
\cite{FalkNeilan}. Other approaches are based on a consistent modification of
the discrete variational formulation and we refer to the overviews
\cite{BoffiBrezziFortin}, \cite{Brenner_Crouzeix}, \cite{ErnGuermondII} for
detailed expositions. In any case, a \textquotedblleft good\textquotedblright%
\ Stokes discretization should have the following features: a) discrete
stability in the form of a discrete inf-sup condition, b) the divergence of
the discrete velocity is zero or very small, c) the Stokes element $\left(
\mathbf{S},M\right)  $ enjoys good approximation properties for the continuous
solution depending on its regularity, d) the element is simple and easy to be implemented.

The Scott-Vogelius element (see \cite{ScottVogelius}, \cite{vogelius1983right}%
) is a very popular element which is based on an appropriate reduction of the
full pressure space in the intuitive element described above. The element is
inf-sup stable, the discrete velocity is divergence-free, and its
implementation very simple. However it suffers from two drawbacks; a) the
discrete inf-sup constant is not mesh-robust: if some vertex becomes nearly
singular -- a geometric notion which will be recalled in the paper and which
is not related to shape regularity -- the discrete inf-sup constant tends to
zero and the discretization becomes increasingly ill-posed; b) in the presence
of super-critical vertices (defined in \eqref{defsupersingular} below), the
approximation property of the pressure space becomes sub-optimal and affects
the accuracy of the discrete pressure significantly. As a remedy for drawback
(a), mesh refinement strategies are proposed in the literature (see \cite[Rem.
2]{Ainsworth_parker_hp_version}) or, alternatively, a very simple modification
of the Scott-Vogelius element which circumvents mesh refinement is introduced
in \cite[Def. 3]{Sauter_eta_wired} 
and called the \textit{pressure-wired} Stokes element. For both methods
the discrete inf-sup constant becomes mesh-robust while the divergence of the
discrete velocity for the pressure-wired Stokes element is not exactly zero
but small (without any regularity requirement on the solution) and controlled
by a parameter $\eta$.

Drawback (b) affects both, the Scott-Vogelius and the pressure-wired Stokes
element: if the mesh contains super-critical vertices, a geometric notion which will be introduced in this paper, the discrete pressure
converges only at a very sub-optimal rate. In this paper, we propose a method
to modify the pressure space so that the pressure converges with optimal rate.
The method is very simple: for both, the Scott-Vogelius and the pressure-wired
Stokes element, the reduction of the\ full pressure space in the most intuitive
Stokes element is formulated via a linear (local) constraint at each critical
vertex written in the form $A_{\mathcal{T},\mathbf{z}}\left(
q\right)  =0$ for a discrete pressure function $q$ and a critical point
$\mathbf{z}$ in the mesh $\mathcal{T}$. This condition is replaced at
super-critical vertices by another linear side condition which we introduce in
this paper. Its definition 
relies on the explicit knowledge of a local basis (set of \emph{critical functions}) for the
orthogonal complement of the reduced pressure space in the full pressure
space. We define a linear injection of the reduced pressure space into the
full pressure space by adding a linear combination of those critical functions
related to super-critical vertices. The coefficients in this
linear combination are given by local linear functionals applied to a pressure
in the reduced pressure space. We emphasize that algorithmically this
modification is simply realized by replacing the
linear constraint $A_{\mathcal{T},\mathbf{z}}\left(  q\right)  =0$ at super-critical vertices with another
local linear constraint and hence the computational complexity and the
dimension of the pressure space stay unchanged. It turns out that the inf-sup stability and the smallness of the divergence of the
discrete velocity is inherited by the resulting element with this modified
pressure space. In addition, the modified pressure space satisfies optimal
approximation properties.

The paper is organized as follows. We introduce the Stokes problem, its
variational form, and the Galerkin discretization in Section
\ref{Sec:StokesP and discretisation}. Section \ref{Sec:Recover approx SV} is
devoted to the Scott-Vogelius element with the modified pressure space. We
introduce the linear functional which serves as the side constraint for the
reduction of the intuitive pressure space, prove that the inf-sup stability of
the Scott-Vogelius element is inherited, and the divergence of the discrete
velocity remains zero. This modification is particularly simple for the
Scott-Vogelius element since it can be realized as a postprocessing step
applied to the original Scott-Vogelius solution. In Subsection
\ref{ApproxSVpressure} and \ref{ApproxSVpressureII} we prove that this
postprocessed pressure converges at optimal rate. In Section \ref{SecPWStokes}
we introduce the pressure-improvement strategy for the pressure-wired Stokes
element and define the modified pressure space. In contrast to the
Scott-Vogelius element, the discretization of the Stokes equation employs the modified pressure space and directly yields the final discrete
solution. We prove in Section \ref{SecStabmodPWStokes} that the inf-sup stability
of the original pressure-wired Stokes element is inherited to its modified
version. In Section \ref{SecPWStokesConv} we will show that our
pressure-improvement strategy applied to the pressure-wired Stokes element
leads to a pressure space with optimal approximation properties. This result
is applied to the modified element in Section \ref{ApproxSVpressureII} and
optimal convergence rates for the discrete solution are shown. It remains to
investigate the divergence of the corresponding discrete velocity which is
considered in Section \ref{DivControl}. The key role is played by the
derivation of an explicit basis representation for the orthogonal complement of
the modified pressure space (Lem. \ref{Prop:Orhtogonal complement Mmod}). This
is used in Theorem \ref{Thm:Divergence estimate} to prove that the smallness of the
velocity divergence is controlled by the parameter $\eta$. For this result no
regularity assumption on the exact solution is imposed. However, if the
continuous solution has some regularity we derive additional convergence rates
with respect to $h$ and $k$ for the smallness of the divergence.

\section{The Stokes problem and its numerical
discretization\label{Sec:StokesP and discretisation}}

Let $\Omega\subset\mathbb{R}^{2}$ be a bounded Lipschitz domain with polygonal
boundary $\partial\Omega$. We consider the numerical solution of the Stokes
equation%
\[%
\begin{array}
[c]{lll}%
-\Delta\mathbf{u}+\nabla p & =\mathbf{f} & \text{in }\Omega,\\
\operatorname*{div}\mathbf{u} & =0 & \text{in }\Omega
\end{array}
\]
with homogeneous Dirichlet boundary conditions for the velocity and the usual
normalization condition for the pressure, namely
\[
\mathbf{u}=\mathbf{0}\quad\text{on }\partial\Omega\quad\text{and\quad}%
\int_{\Omega}p=0.
\]
Throughout this paper, standard notation applies for real-valued Lebesgue and
Sobolev spaces. Let $H_{0}^{1}\left(  \Omega\right)  $ be the closure of the
space of infinitely smooth, compactly supported functions with respect
to the $H^{1}\left(  \Omega\right)  $ norm. Its dual space is given by
$H^{-1}\left(  \Omega\right)  :=H_{0}^{1}\left(  \Omega\right)  ^{\prime}$.
The scalar product and norm in $L^{2}\left(  \Omega\right)  $ are written as%
\[
\left(  u,v\right)  _{L^{2}\left(  \Omega\right)  }:=\int_{\Omega}%
uv\quad\text{and}\quad\left\Vert u\right\Vert _{L^{2}\left(  \Omega\right)
}:=\left(  u,u\right)  _{L^{2}\left(  \Omega\right)  }^{1/2}.
\]
Vector-valued and $2\times2$ tensor-valued analogues of these function spaces
are denoted by bold and blackboard bold letters, e.g., $\mathbf{H}^{s}\left(
\Omega\right)  =\left(  H^{s}\left(  \Omega\right)  \right)  ^{2}$ and
$\mathbb{H}^{s}\left(  \Omega\right)  =\left(  H^{s}\left(  \Omega\right)
\right)  ^{2\times2}$ and analogously for other quantities.

The $\mathbf{L}^{2}\left(  \Omega\right)  $ scalar product and norm for
vector-valued functions are given by%
\[
\left(  \mathbf{u},\mathbf{v}\right)  _{\mathbf{L}^{2}\left(  \Omega\right)
}:=\int_{\Omega}\left\langle \mathbf{u},\mathbf{v}\right\rangle \quad
\text{and\quad}\left\Vert \mathbf{u}\right\Vert _{\mathbf{L}^{2}\left(
\Omega\right)  }:=\left(  \mathbf{u},\mathbf{u}\right)  _{\mathbf{L}%
^{2}\left(  \Omega\right)  }^{1/2},
\]
with the standard Euclidean scalar product $\left\langle \cdot,\cdot\right\rangle$.
In a similar fashion, we define the scalar product and norm in
$\mathbb{L}^{2}\left(  \Omega\right)  $ by%
\[
\left(  \mathbf{G},\mathbf{H}\right)  _{\mathbb{L}^{2}\left(  \Omega\right)
}:=\int_{\Omega}\left\langle \mathbf{G},\mathbf{H}\right\rangle \quad
\text{and\quad}\left\Vert \mathbf{G}\right\Vert _{\mathbb{L}^{2}\left(
\Omega\right)  }:=\left(  \mathbf{G},\mathbf{G}\right)  _{\mathbb{L}%
^{2}\left(  \Omega\right)  }^{1/2} \quad\forall\mathbf{G},\mathbf{H}
\in\mathbb{L}^{2}\left(  \Omega\right)  ,
\]
where $\left\langle \mathbf{G},\mathbf{H}\right\rangle =\sum_{i,j=1}%
^{2}G_{i,j}H_{i,j}$. Finally, let $L_{0}^{2}\left(  \Omega\right)  :=\left\{
u\in L^{2}\left(  \Omega\right)  :\int_{\Omega}u=0\right\}  $. We introduce
the bilinear forms $a:\mathbf{H}^{1}\left(  \Omega\right)  \times
\mathbf{H}^{1}\left(  \Omega\right)  \rightarrow\mathbb{R}$ and $b:\mathbf{H}%
^{1}\left(  \Omega\right)  \times L^{2}\left(  \Omega\right)  \rightarrow
\mathbb{R}$ by%
\begin{equation}%
\begin{array}
[c]{ccc}%
a\left(  \mathbf{u},\mathbf{v}\right)  :=\left(  \nabla\mathbf{u}%
,\nabla\mathbf{v}\right)  _{\mathbb{L}^{2}\left(  \Omega\right)  } & \text{
and } & b\left(  \mathbf{u},p\right)  =\left(  \operatorname{div}%
\mathbf{u},p\right)  _{L^{2}\left(  \Omega\right)  }\label{defabili}%
\end{array}
\end{equation}
with the derivative
$\nabla\mathbf{v}$ and the divergence $\operatorname{div} \mathbf{v}$ of any
$\mathbf{v} \in\mathbf{H}^{1} \left(  \Omega\right)  $. Given a source
$\mathbf{F}\in\mathbf{H}^{-1}\left(  \Omega\right)  $, the variational form of
the stationary Stokes problem seeks $\left(  \mathbf{u},p\right)
\in\mathbf{H}_{0}^{1}\left(  \Omega\right)  \times L_{0}^{2}\left(
\Omega\right)  $ such that
\begin{equation}%
\begin{array}
[c]{lll}%
a\left(  \mathbf{u},\mathbf{v}\right)  -b\left(  \mathbf{v},p\right)  &
=\mathbf{F}\left(  \mathbf{v}\right)  & \forall\mathbf{v}\in\mathbf{H}_{0}%
^{1}\left(  \Omega\right)  ,\\
b\left(  \mathbf{u},q\right)  & =0 & \forall q\in L_{0}^{2}\left(
\Omega\right)  .
\end{array}
\label{varproblemstokes}%
\end{equation}
Concerning the well-posedenss of \eqref{varproblemstokes}, we refer, e.g., to
\cite{Girault86} for details. In this paper, we consider a conforming Galerkin
discretization of (\ref{varproblemstokes}) by a pair $\left(  \mathbf{S}%
,M\right)  $ of finite dimensional subspaces of the continuous solution spaces
$\left(  \mathbf{H}_{0}^{1}\left(  \Omega\right)  ,L_{0}^{2}\left(
\Omega\right)  \right)  $. For any given $\mathbf{F}\in\mathbf{H}^{-1}\left(
\Omega\right)  $ the weak formulation seeks $\left(  \mathbf{u}_{\mathbf{S}%
},p_{M}\right)  \in\mathbf{S}\times M\;$such that%
\begin{equation}%
\begin{array}
[c]{lll}%
a\left(  \mathbf{u}_{\mathbf{S}},\mathbf{v}\right)  -b\left(  \mathbf{v}%
,p_{M}\right)  & =\mathbf{F}\left(  \mathbf{v}\right)  & \forall\mathbf{v}%
\in\mathbf{S},\\
b\left(  \mathbf{u}_{\mathbf{S}},q\right)  & =0 & \forall q\in M.
\end{array}
\label{discrStokes}%
\end{equation}
It is well known that the bilinear form $a\left(  \cdot,\cdot\right)  $ is
symmetric, continuous, and coercive so that problem (\ref{discrStokes}) is
well-posed if the bilinear form $b\left(  \cdot,\cdot\right)  $ satisfies the
inf-sup condition for $\left(  \mathbf{S}, M \right)  $.

\begin{definition}
Let $\mathbf{S}$ and $M$ be finite-dimensional subspaces of $\mathbf{H}%
_{0}^{1}\left(  \Omega\right)  $ and $L_{0}^{2}\left(  \Omega\right)  $. The
pair $\left(  \mathbf{S},M\right)  $ is called \emph{inf-sup stable} if the
\emph{inf-sup constant}
is positive:
\begin{equation}
\beta\left(  \mathbf{S},M\right)  :=\inf_{q\in M\backslash\left\{  0\right\}
}\sup_{\mathbf{v}\in\mathbf{S}\backslash\left\{  \mathbf{0}\right\}  }%
\frac{\left(  q,\operatorname*{div}\mathbf{v}\right)  _{L^{2}\left(
\Omega\right)  }}{\left\Vert \mathbf{v}\right\Vert _{\mathbf{H}^{1}\left(
\Omega\right)  }\left\Vert q\right\Vert _{L^{2}\left(  \Omega\right)  }}>0.
\label{infsupcond}%
\end{equation}

\end{definition}

\section{Pressure improved Scott-Vogelius element
\label{Sec:Recover approx SV}}

Let $\mathcal{T}$ be a conforming, shape-regular triangulation of the domain
$\Omega$ into closed triangles $K\in\mathcal{T}$ with diameter $h_{K}$. The
set of vertices is given by $\mathcal{V}\left(  \mathcal{T}\right)  $ and the
additional subscripts $\mathcal{V}_{\Omega}\left(  \mathcal{T}\right)  $,
$\mathcal{V}_{\partial\Omega}\left(  \mathcal{T}\right)  $ specify whether a
vertex is located in the domain or on its boundary. The set of edges is
denoted by $\mathcal{E}\left(  \mathcal{T}\right)  $ and the same subscript
convention as for the vertices applies. For a vertex $\mathbf{z}\in
\mathcal{V}\left(  \mathcal{T}\right)  $, the local vertex patch is given by
\begin{equation}%
\begin{array}
[c]{ll}%
\mathcal{T}_{\mathbf{z}}:=\left\{  K\in\mathcal{T}\mid\mathbf{z}\in K\right\}
& \text{and\quad}\omega_{\mathbf{z}}:=%
{\displaystyle\bigcup\limits_{K\in\mathcal{T}_{\mathbf{z}}}}
K
\end{array}
\label{nodalpatch}%
\end{equation}
with the local mesh width $h_{\mathbf{z}}:=\max\left\{  h_{K}:K\in
\mathcal{T}_{\mathbf{z}}\right\}  $. For any vertex $\mathbf{z}\in
\mathcal{V}\left(  \mathcal{T}\right)  $, we fix a local counterclockwise
numbering of the $N_{\mathbf{z}}:=\operatorname{card}\mathcal{T}_{\mathbf{z}}$
triangles in
\begin{equation}
\mathcal{T}_{\mathbf{z}}=\left\{  K_{j}:1\leq j\leq N_{\mathbf{z}}\right\}  .
\label{eqn:enum_Tz}%
\end{equation}
A triangle neighborhood of some triangle $K\in\mathcal{T}$ is given by%
\begin{equation}
\omega\left(  K\right)  :=%
{\displaystyle\bigcup\limits_{\substack{K^{\prime}\in\mathcal{T}\\K^{\prime
}\cap K\neq\emptyset}}}
K^{\prime}. \label{defomegaK}%
\end{equation}
The shape-regularity constant
\begin{equation}
\gamma_{\mathcal{T}}:=\max_{K\in\mathcal{T}}\frac{h_{K}}{\rho_{K}}
\label{defgammat}%
\end{equation}
relates the local mesh width $h_{K}$ with the diameter $\rho_{K}$ of the
largest inscribed ball in an element $K\in\mathcal{T}$. The global mesh width
is given by $h_{\mathcal{T}}:=\max\left\{  h_{K}:K\in\mathcal{T}\right\}  $.
For a subset $M\subset\mathbb{R}^{2}$, we denote the area of $M$ by
$\left\vert M\right\vert $. Let $\mathbb{P}_{k}(K)$ denote the space of
polynomials on $K\in\mathcal{T}$ with total degree smaller than or equal to
$k\in\mathbb{N}_{0}$ and define%
\begin{equation}%
\begin{array}
[c]{l}%
\mathbb{P}_{k}\left(  \mathcal{T}\right)  :=\left\{  q\in L^{2}\left(
\Omega\right)  \mid\forall K\in\mathcal{T}:\left.  q\right\vert _{\overset
{\circ}{{K}}}\in\mathbb{P}_{k}\left(  \overset{\circ}{{K}}\right)  \right\}
,\\
\mathbb{P}_{k,0}\left(  \mathcal{T}\right)  :=\mathbb{P}_{k}\left(
\mathcal{T}\right)  \cap L_{0}^{2}\left(  \Omega\right)  =\left\{
q\in\mathbb{P}_{k}\left(  \mathcal{T}\right)  \mid\int_{\Omega}q=0\right\}
,\\
S_{k}\left(  \mathcal{T}\right)  :=\mathbb{P}_{k}\left(  \mathcal{T}\right)
\cap H^{1}(\Omega),\\
S_{k,0}\left(  \mathcal{T}\right)  :=S_{k}\left(  \mathcal{T}\right)  \cap
H_{0}^{1}\left(  \Omega\right)  .
\end{array}
\label{Pkdefs}%
\end{equation}
For any $q\in\mathbb{P}_{k}\left(  \mathcal{T}\right)  $, we write
\begin{equation}
q_{\operatorname*{mvz}}:=q-\overline{q}\quad\text{with the integral mean\quad
}\overline{q}:=\frac{1}{\left\vert \Omega\right\vert }\int_{\Omega}q.
\label{defq0qbar}%
\end{equation}

\subsection{The Scott-Vogelius element}

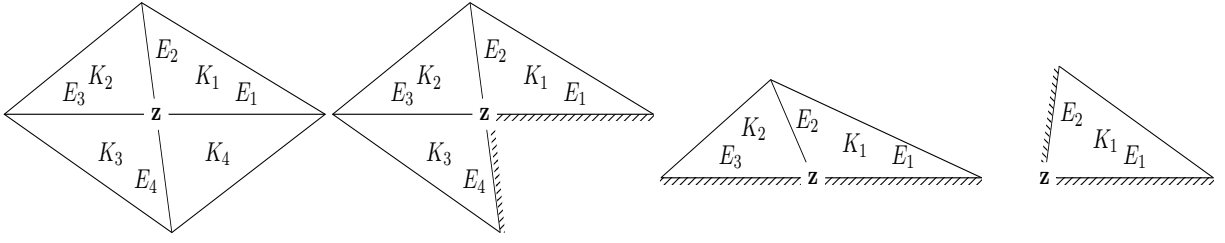
\begin{figure}[ptb]
\center
\resizebox{\textwidth}{8em}{
	\begin{tikzpicture}[scale=0.8]
				\draw (0,0) node[] (z) {$\mathbf{z}$};
		\coordinate (A) at (4,0);
		\coordinate (B) at ($(z)!0.5!100:(A)$);
		\coordinate (C) at ($(z)!0.9!180:(A)$);
		\coordinate (D) at ($(z)!0.53!280:(A)$);
		\draw (A) -- (B) -- (C) -- (D) -- cycle;
		\draw (z) -- (A) node[midway,above] {$E_1$};
		\draw (z) -- (B) node[midway,right] {$E_2$};
		\draw (z) -- (C) node[midway,above] {$E_3$};
		\draw (z) -- (D) node[midway,left] {$E_4$};
		\node (T1) at (barycentric cs:z=1,A=1,B=1) {$K_1$};
		\node (T1) at (barycentric cs:z=1,C=1,B=1) {$K_2$};
		\node (T1) at (barycentric cs:z=1,A=1,D=1) {$K_4$};
		\node (T1) at (barycentric cs:z=1,C=1,D=1) {$K_3$};
	\end{tikzpicture}
	\begin{tikzpicture}[scale=0.8]
				\draw (0,0) node[] (z) {$\mathbf{z}$};
		\coordinate (A) at (4,0);
		\coordinate (B) at ($(z)!0.5!100:(A)$);
		\coordinate (C) at ($(z)!0.9!180:(A)$);
		\coordinate (D) at ($(z)!0.53!280:(A)$);
		\draw (A) -- (B) -- (C) -- (D);
		\draw (z) -- (A) node[midway,above] {$E_1$};
		\draw (z) -- (B) node[midway,right] {$E_2$};
		\draw (z) -- (C) node[midway,above] {$E_3$};
		\draw (z) -- (D) node[midway,left] {$E_4$};
		\node (T1) at (barycentric cs:z=1,A=1,B=1) {$K_1$};
		\node (T1) at (barycentric cs:z=1,C=1,B=1) {$K_2$};
				\node (T1) at (barycentric cs:z=1,C=1,D=1) {$K_3$};
		\fill[pattern=north east lines, ] ([xshift=0.35em]z.south) -- (D) -- ([xshift=.3em]D) -- ([xshift=.65em]z.south);
		\fill[pattern=north east lines, ] (z.east) rectangle ([yshift=-.3em]A);
	\end{tikzpicture}
	\begin{tikzpicture}[scale=0.8]
				\draw (0,0) node[] (z) {$\mathbf{z}$};
		\coordinate (A) at (4,0);
		\coordinate (B) at ($(z)!0.5!120:(A)$);
		\coordinate (C) at ($(z)!0.9!180:(A)$);
		\coordinate (D) at ($(z)!0.53!300:(A)$);
		\draw (A) -- (B) -- (C);		\draw (z) -- (A) node[midway,above] {$E_1$};
		\draw (z) -- (B) node[midway,right] {$E_2$};
		\draw (z) -- (C) node[midway,above] {$E_3$}; 
		\node (T1) at (barycentric cs:z=1,A=1,B=1) {$K_1$};
		\node (T1) at (barycentric cs:z=1,C=1,B=2) {$K_2$};
		\phantom{		\node (T1) at (barycentric cs:z=1,A=1,D=1) {$K_4$};
		\node (T1) at (barycentric cs:z=1,C=1,D=1) {$K_3$};}
		\fill[pattern=north east lines, ] (z.east) rectangle ([yshift=-.3em]A);
		\fill[pattern=north east lines, ] (z.west) rectangle ([yshift=-.3em]C);
	\end{tikzpicture}
	\begin{tikzpicture}[scale=0.8]
				\draw (0,0) node[] (z) {$\mathbf{z}$};
		\coordinate (A) at (4,0);
		\coordinate (B) at ($(z)!0.5!80:(A)$);
		\coordinate (C) at ($(z)!0.9!180:(A)$);
		\coordinate (D) at ($(z)!0.53!300:(A)$);
		\draw (A) -- (B);		\draw (z) -- (A) node[midway,above] {$E_1$};
		\draw (z) -- (B) node[midway,right] {$E_2$};
		\phantom{		\node (T1) at (barycentric cs:z=1,A=1,D=1) {$K_4$};
		\node (T1) at (barycentric cs:z=1,C=1,D=1) {$K_3$};}
		\node (T1) at (barycentric cs:z=1,A=1,B=1) {$K_1$};
		\fill[pattern=north east lines, ] (z.east) rectangle ([yshift=-.3em]A);
		\fill[pattern=north west lines, ] ([xshift=0.08em]z.north) -- (B) -- ([xshift=-.3em]B) -- ([xshift=-.22em]z.north);
\end{tikzpicture}}
\caption{Vertex patch for an interior singular vertex
$\mathbf{z}\in\mathcal{V}_{\Omega}(\mathcal{T})$ with $N_{\mathbf{z}}=4$
(resp.~boundery singular vertex $\mathbf{z}\in\mathcal{V}_{\partial\Omega
}(\mathcal{T})$ with $N_{\mathbf{z}}=1,2,3$) triangles}%
\label{fig:Singular vertex_patch}%
\end{figure}

It is well known that the most intuitive Stokes element $\left(
\mathbf{S}_{k,0}\left(  \mathcal{T}\right)  ,\mathbb{P}_{k-1,0}\left(
\mathcal{T}\right)  \right)  $ is in general unstable. The analysis in
 \cite{ScottVogelius}, \cite{vogelius1983right} for $k\geq4$ relates the
instability of $\left(  \mathbf{S}_{k,0}\left(  \mathcal{T}\right)
,\mathbb{P}_{k-1,0}\left(  \mathcal{T}\right)  \right)  $ to the presence of
\emph{critical} or \emph{singular points} in the mesh $\mathcal{T}$. The
following definition is illustrated by Fig. \ref{fig:Singular vertex_patch}.

\begin{definition}
\label{DefCritpoint} The local measure of singularity $\Theta\left(
\mathbf{z}\right)  $ at $\mathbf{z}\in\mathcal{V}(\mathcal{T})$ reads
\begin{equation}
\Theta\left(  \mathbf{z}\right)  :=%
\begin{cases}
\max\left\{  \left.  \left\vert \sin\left(  \theta_{i}+\theta_{i+1}\right)
\right\vert \;\right\vert \;0\leq i\leq N_{\mathbf{z}}\right\}  & \text{if
}\mathbf{z}\in\mathcal{V}_{\Omega}\left(  \mathcal{T}\right)  ,\\
\max\left\{  \left.  \left\vert \sin\left(  \theta_{i}+\theta_{i+1}\right)
\right\vert \;\right\vert \;0\leq i\leq N_{\mathbf{z}}-1\right\}  & \text{if
}\mathbf{z}\in\mathcal{V}_{\partial\Omega}(\mathcal{T})\wedge N_{\mathbf{z}%
}>1,\\
0 & \text{if }\mathbf{z\in}\mathcal{V}_{\partial\Omega}(\mathcal{T})\wedge
N_{\mathbf{z}}=1,
\end{cases}
\label{defthetaz}%
\end{equation}
where the angles $\theta_{j}$ in $K_{j}\in\mathcal{T}%
_{\mathbf{z}}$ at $\mathbf{z}$ are numbered counterclockwise
from $1\leq j\leq N_{\mathbf{z}}$ (see (\ref{eqn:enum_Tz})) and cyclic
numbering is applied: $\theta_{N_{\mathbf{z}}+1}:=\theta_{1}$ if the patch is
closed, i.e., $\mathbf{z}\in\mathcal{V}_{\Omega}\left(  \mathcal{T}\right)  $.
A vertex $\mathbf{z}\in\mathcal{V}\left(  \mathcal{T}\right)  $ with
$\Theta\left(  \mathbf{z}\right)  =0$ is a \emph{singular vertex} and the set
of all singular vertices is%
\[
\mathcal{C}_{\mathcal{T}}:=\left\{  \mathbf{z}\in\mathcal{V}\left(
\mathcal{T}\right)  \mid\Theta\left(  \mathbf{z}\right)  =0\right\}  .
\]
The global \emph{measure of singularity} of the mesh $\mathcal{T}$ is
\begin{equation}
\Theta_{\min}:=\min_{\mathbf{z}\in\mathcal{V}\left(  \mathcal{T}\right)
\backslash\mathcal{C}_{\mathcal{T}}}\Theta\left(  \mathbf{z}\right)  .
\label{thetamin}%
\end{equation}
For any vertex $\mathbf{z}\in\mathcal{V}\left(  \mathcal{T}\right)  $ and all
$q\in\mathbb{P}_{k-1}\left(  \mathcal{T}\right)  $, the functional
$A_{\mathcal{T},\mathbf{z}}$ is the alternating sum
\begin{equation}
A_{\mathcal{T},\mathbf{z}}\left(  q\right)  :=\sum_{\ell=1}^{N_{\mathbf{z}}%
}\left(  -1\right)  ^{\ell}\left(  \left.  q\right\vert _{K_{\ell}}\right)
\left(  \mathbf{z}\right)   \label{DefATz}%
\end{equation}
over the triangles $K_{\ell}\in\mathcal{T}_{\mathbf{z}}$ for $1\leq\ell\leq
N_{\mathbf{z}}$.
\end{definition}

In \cite{ScottVogelius} and \cite{vogelius1983right}, the space%
\begin{equation}
M_{0,k-1}\left(  \mathcal{T}\right)  :=\left\{  \left.  q\in\mathbb{P}%
_{k-1,0}\left(  \mathcal{T}\right)  \;\right\vert \;\forall\mathbf{z}%
\in\mathcal{C}_{\mathcal{T}}:\ A_{\mathcal{T},\mathbf{z}}\left(  q\right)
=0\right\}  \label{SVpressure}%
\end{equation}
was introduced and used in the definition of the Scott-Vogelius element
$\left(  \mathbf{S}_{k,0}\left(  \mathcal{T}\right)  ,M_{0,k-1}\left(
\mathcal{T}\right)  \right)  $. It was proven in \cite{vogelius1983right} that
this element enjoys two important properties of a \textquotedblleft
good\textquotedblright\ Stokes element: a) inf-sup stability which follows
from%
\begin{equation}
\operatorname{div}\mathbf{S}_{k,0}\left(  \mathcal{T}\right)  =M_{0,k-1}%
\left(  \mathcal{T}\right)  \label{Eq:Def SV M}%
\end{equation}
and b) the discrete velocity is divergence free. A third important property
certainly is the approximation property of the discrete spaces. Since
$\mathbf{S}_{k,0}\left(  \mathcal{T}\right)  $ is a standard finite element
space its approximation property is well known. For the pressure space, however,
the approximation property might deteriorate in the vicinity of vertices
having a particular type of vertex patch $\mathcal{T}_{\mathbf{z}}$ as
explained next. We start with an observation for $q\in M_{0,k-1}\left(
\mathcal{T}\right)  $. Suppose $\mathbf{z}\in\mathcal{C}_{\mathcal{T}}%
\cap\partial\Omega$ is a singular boundary vertex with an odd number of
neighboring triangles (type 2 or 4 in Figure~\ref{fig:Singular vertex_patch}).
The side condition in the definition of the pressure space $A_{\mathcal{T}%
,\mathbf{z}}\left(  q\right)  =0$ reveals the following
implication
\begin{equation}
q\in M_{0,k-1}\left(  \mathcal{T}\right)  \text{ is continuous at }%
\mathbf{z}\implies q\left(  \mathbf{z}\right)  =0. \label{qz0}%
\end{equation}
Since the exact pressure does not vanish at theses points in
general, we cannot expect a good approximation property of $M_{0,k-1}\left(
\mathcal{T}\right)  $ in neighborhoods of such vertices. In particular, the
$L^{\infty}$ norm $\left\Vert p-p_{M}\right\Vert _{L^{\infty}\left(
\Omega\right)  }$ of a smooth pressure $p\in C^{0}\left(  \overline{\Omega
}\right)  $ is at least $\left\vert p\left(  \mathbf{z}\right)  \right\vert $,
independent of the triangulation $\mathcal{T}$. This is a drawback and we will
present a strategy to modify the pressure space in these vertices such that
standard approximation properties hold. The singular vertices responsible for
the deficiency in the approximation properties of $M_{0,k-1}\left(
\mathcal{T}\right)  $ are collected in the set of \emph{super-critical}
vertices:%
\begin{equation}
\mathcal{SC}_{\mathcal{T}}:=\left\{  \left.  \mathbf{z}\in\mathcal{C}%
_{\mathcal{T}}\;\right\vert \;N_{\mathbf{z}}\text{ is odd}\right\}  .
\label{defsupersingular}%
\end{equation}

\begin{remark}
\label{Rem:even singular vertex} It was shown in \cite{Ainsworth_parker_II}
that the space
\begin{equation}
\widetilde{Q}_{0}^{h,k}:=\left\{  \left.  q\in\mathbb{P}_{k-1,0}\left(
\mathcal{T}\right)  \;\right\vert \;q\text{ is }C^{0}\text{ at all
vertices}\right\}  \label{defQtilde}%
\end{equation}
has optimal approximation properties. Given any $q\in\widetilde{Q}_{0}^{h,k}$,
Definition \ref{DefCritpoint} provides $A_{\mathcal{T},\mathbf{z}}\left(
q\right)  =0 $ for all $\mathbf{z}\in\mathcal{C}_{\mathcal{T}}$ with
$N_{\mathbf{z}}$ even. Therefore, we only expect the loss of approximability in
the vertices in $\mathcal{SC}_{\mathcal{T}}$.
\end{remark}

As discussed in the introduction, we will present a simple modification
strategy to remedy the reduced approximation properties. For the
Scott-Vogelius element, this strategy leads to a local postprocessing
step and recovers the optimal convergence order for the pressure approximation.

\subsection{A simple postprocessing strategy}

Since the discrete pressure space
$M_{0,k-1}\left(  \mathcal{T}\right)  $ is defined by restricting the full
pressure space $\mathbb{P}_{k-1,0}\left(  \mathcal{T}\right)  $ at singular
vertices via the functional $A_{\mathcal{T},\mathbf{z}}$, the complement of
$M_{0,k-1}\left(  \mathcal{T}\right)  $ in $\mathbb{P}_{k-1,0}\left(
\mathcal{T}\right)  $ is non-trivial in the presence of singular vertices and
has been described in \cite{Baran_Stoyan} and \cite[Def. 3.11]{CCSS_CR_1}.

\begin{definition}
\label{Def:bz} For a vertex $\mathbf{z}\in\mathcal{V}\left(  \mathcal{T}%
\right)  $, the \emph{critical function }$b_{k-1,\mathbf{z}}\in\mathbb{P}%
_{k-1}\left(  \mathcal{T}\right)  $ is given by%
\begin{equation}
b_{k-1,\mathbf{z}}=\sum_{\ell=1}^{N_{\mathbf{z}}}\frac{\left(  -1\right)
^{k-1+\ell}}{\left\vert K_{\ell}\right\vert }P_{k-1}^{\left(  0,2\right)
}\left(  1-2\lambda_{K_{\ell},\mathbf{z}}\right)  \chi_{K_{\ell}%
},\label{Eq:Def bz}%
\end{equation}
where $\chi_{K_{\ell}}$ is the characteristic function of the triangle
$K_{\ell}\in\mathcal{T}_{\mathbf{z}}$ and $P_{k-1}^{\left(  \alpha
,\beta\right)  }$ is the Jacobi polynomial used here for the parameters
$\alpha=0$, $\beta=2$ (cf., \cite[Table 18.3.1]{NIST:DLMF}). The barycentric
coordinate for a triangle $K\in\mathcal{T}\left(  \mathbf{z}\right)  $
corresponding to the vertex $\mathbf{z}$ is denoted by $\lambda_{K,\mathbf{z}%
}$.
\end{definition}

The following properties of the critical functions $b_{k-1,\mathbf{z}}$ were
proven in \cite[Lem. 6]{Sauter_eta_wired}.

\begin{proposition}
\label{Propcritfct}The critical functions $b_{k-1,\mathbf{z}}$ satisfy for all
$K_{\ell}\in\mathcal{T}_{\mathbf{z}}:$%
\begin{subequations}
\begin{align}
\left(  b_{k-1,\mathbf{z}},1\right)  _{L^{2}\left(  K_{\ell}\right)  }  &
=\left(  -1\right)  ^{\ell}\binom{k+1}{2}^{-1},\label{bkintegral}\\
\left\Vert b_{k-1,\mathbf{z}}\right\Vert _{L^{2}\left(  K_{\ell}\right)
}^{2}  &  =\left\vert K_{\ell}\right\vert ^{-1},\label{bkL2}\\
\left.  b_{k-1,\mathbf{z}}\right\vert _{K_{\ell}}\left(  \mathbf{y}\right)
&  =\frac{\left(  -1\right)  ^{\ell}}{\left\vert K_{\ell}\right\vert }\left\{
\begin{array}
[c]{ll}%
\binom{k+1}{2} & \text{if }\mathbf{y}=\mathbf{z},\\
\left(  -1\right)  ^{k-1} & \text{otherwise,}%
\end{array}
\right. \label{bkpointv}\\
A_{\mathcal{T},\mathbf{z}}\left(  b_{k-1,\mathbf{z}}\right)   &  =\binom
{k+1}{2}\left\Vert b_{k-1,\mathbf{z}}\right\Vert _{L^{2}\left(  \omega
_{\mathbf{z}}\right)  }^{2}. \label{Azpropb}%
\end{align}
For $k\geq2$, the set $\left\{  b_{k-1,\mathbf{z}}:\mathbf{z}\in
\mathcal{V}\left(  \mathcal{T}\right)  \right\}  \cup\left\{  1\right\}  $ is
linearly independent and%
\end{subequations}
\begin{equation}
\frac{3}{4}\left\Vert \sum_{\mathbf{z}\in\mathcal{V}\left(  K\right)
}c_{\mathbf{z}}b_{k,\mathbf{z}}\right\Vert _{L^{2}\left(  K\right)  }^{2}%
\leq\left\vert K\right\vert ^{-1}\sum_{\mathbf{z}\in\mathcal{V}\left(
K\right)  }c_{\mathbf{z}}^{2}\leq\frac{12}{7}\min_{\alpha\in\mathbb{R}%
}\left\Vert \sum_{\mathbf{z}\in\mathcal{V}\left(  K\right)  }c_{\mathbf{z}%
}b_{k,\mathbf{z}}-\alpha\right\Vert _{L^{2}\left(  K\right)  }^{2}.
\label{stabest}%
\end{equation}

\end{proposition}

The functions $b_{k-1,\mathbf{z}}$ characterise the
orthogonal complement of $M_{0,k-1}\left(  \mathcal{T}\right)$ in
$\mathbb{P}_{k-1,0}\left(  \mathcal{T}\right)$.
Let $\oplus$ denote the direct sum of vector spaces and recall $\left(
\cdot\right)  _{\operatorname*{mvz}}$ from~\eqref{defq0qbar}.

\begin{proposition}[{\cite[Lem. 3.13]{CCSS_CR_1}}]
\label{PropOrthCompl}

The decomposition%
\[
\mathbb{P}_{k-1,0}\left(  \mathcal{T}\right)  =M_{0,k-1}\left(  \mathcal{T}%
\right)  \oplus\operatorname*{span}\left\{  \left(  b_{k-1,\mathbf{z}}\right)
_{\operatorname*{mvz}}\mid\mathbf{z}\in\mathcal{C}_{\mathcal{T}}\right\} 
\]
is $L^2$ orthogonal, i.e., any $q_M\in M_{0,k-1}(\mathcal{T})$ satisfies $(q_M, \left(  b_{k-1,\mathbf{z}}\right)
_{\operatorname*{mvz}})_{L^2(\Omega)}=0$ for all $\mathbf{z}\in \mathcal{C}_{\mathcal{T}}$.
\end{proposition}
\begin{proof}
	This is a direct consequence of~\cite[Lem. 3.13]{CCSS_CR_1}, the definition of 
	$\left( \cdot\right)  _{\operatorname*{mvz}} $, and the integral mean zero condition in $M_{0,k-1}(\mathcal{T})$; further
	details are omitted.
\end{proof}

We will employ a continuous, linear functional
\begin{equation}
f_{\mathbf{z}}:\mathbb{P}_{k-1}\left(  \mathcal{T}\right)  \rightarrow
\mathbb{R\quad}\text{with }\left\vert f_{\mathbf{z}}\left(  q\right)
\right\vert \leq\frac{C_{f_{\mathbf{z}}}}{\left\Vert b_{k-1,\mathbf{z}%
}\right\Vert _{L^{2}\left(  \Omega\right)  }}\left\Vert q\right\Vert
_{L^{2}\left(  \Omega\right)  } \label{estfz}%
\end{equation}
to define the modified pressure space $M_{0,k-1}^{\operatorname{mod}}\left(
\mathcal{T}\right)  \subset\mathbb{P}_{k-1,0}\left(  \mathcal{T}\right)  $ by:%
\begin{equation}
M_{0,k-1}^{\operatorname{mod}}\left(  \mathcal{T}\right)  :=\left\{
q+\sum_{\mathbf{z}\in\mathcal{SC}_{\mathcal{T}}}f_{\mathbf{z}}\left(
q\right)  \left(  b_{k-1,\mathbf{z}}\right)  _{\operatorname*{mvz}}\;:\;q\in
M_{0,k-1}\left(  \mathcal{T}\right)  \right\}  . \label{Eq:Def SV Mmod}%
\end{equation}
The modification in $M_{0,k-1}^{\operatorname{mod}}\left(  \mathcal{T}\right)$ for general $f_{\mathbf{z}}$ overcomes
the implication~\eqref{qz0} that results in the suboptimal approximation properties of $M_{0,k-1}\left(
\mathcal{T}\right)$ in the presence of super critical vertices and defines a novel discretisation
of~\eqref{varproblemstokes}.
A \emph{good} choice of $f_{\mathbf{z}}$ discussed in Subsection~\ref{sub:assumptions_on_fz} even enables optimal approximation properties of
$M_{0,k-1}^{\operatorname{mod}}\left(  \mathcal{T}\right)$ in the sense of Theorem~\ref{Thm:Recoverd bestapprox} below.

\begin{definition}
	Given $k\in\mathbb N$ and functionals $f_{\mathbf{z}}$ with (\ref{estfz}) for all $\mathbf{z}\in\mathcal{C}_{\mathcal{T}}$ and
	$M_{0,k-1}^{\operatorname{mod}}(\mathcal{T})$ from~\eqref{Eq:Def SV Mmod}, the
pressure-modified Scott-Vogelius element is given by the pair $\left(
\mathbf{S}_{k,0}\left(  \mathcal{T}\right)  ,M_{0,k-1}^{\operatorname{mod}%
}\left(  \mathcal{T}\right)  \right)  $.
\end{definition}

An important point is that the pressure-modified Scott-Vogelius solution is a simple post-processing of the classical Scott-Vogelius solution
with possibly better approximation properties.
\begin{theorem}
\label{Thm:ApproxPropertie pp-pressure}
Let $k\geq4$ and $\mathbf{F}\in\mathbf{H}^{-1}\left(
\Omega\right)$.
Then $\left(  \mathbf{u}_{\mathbf{S}},p_{M}\right)\in \mathbf{S}_{k,0}\left(  \mathcal{T}\right)
\times M_{0,k-1}\left(  \mathcal{T}\right)$ solves~\eqref{discrStokes} in
$(\mathbf{S}_{k,0}\left(  \mathcal{T}\right),M_{0,k-1}\left(  \mathcal{T}\right))$
 if and only if  
$\left(  \mathbf{u}_{\mathbf{S}},p_{M}^*\right)$ solves~\eqref{discrStokes} in
$(\mathbf{S}_{k,0}\left(  \mathcal{T}\right),M_{0,k-1}^{\operatorname{mod}}\left(  \mathcal{T}\right))$
with
\begin{equation}
p_{M}^{\ast}:=p_{M}+\sum_{\mathbf{z}\in\mathcal{SC}_{\mathcal{T}}%
}f_{\mathbf{z}}\left(  p_{M}\right)  \left(  b_{k-1,\mathbf{z}}\right)
_{\operatorname*{mvz}}\in M_{0,k-1}^{\operatorname{mod}}\left(  \mathcal{T}\right).\label{defpostprocessedpressure}%
\end{equation}
Given $C_f\coloneqq 1+\sum^{}_{\mathbf{z}\in\mathcal{S}\mathcal{C}_{\mathcal{T}}} C_{f_{\mathbf{z}}}$, the
modified discrete pressure $p_{M}^{\ast}$ satisfies the error estimate
\begin{align}
\left\Vert p-p_{M}^{\ast}\right\Vert _{L^{2}\left(  \Omega\right)  }\leq &
\frac{C_{\operatorname{vel}}C_f^{2}}{\Theta_{\min}^{2}}%
\inf_{\mathbf{v}\in\mathbf{S}_{k,0}\left(  \mathcal{T}\right)  }\left\Vert
\mathbf{u}-\mathbf{v}\right\Vert _{H^{1}\left(  \Omega\right)  }
  +\frac{C_{\operatorname{pres}}C_f}{\Theta_{\min}}\inf_{q\in M_{0,k-1}%
^{\operatorname{mod}}\left(  \mathcal{T}\right)  }\left\Vert p-q\right\Vert
_{L^{2}\left(  \Omega\right)  }\label{Eq:QuasiBestApprox of Mmod}%
\end{align}
with the continuous solution
$\left(  \mathbf{u},p\right)  \in\mathbf{H}_{0}^{1}\left(  \Omega\right)
\times
L_{0}^{2}\left(  \Omega\right)  $
to~\eqref{varproblemstokes}.
The positive constants
$C_{\operatorname{vel}},C_{\operatorname{pres}}$ only depend on the
shape-regularity of the mesh and the domain $\Omega$.
\end{theorem}
%
%

The proof of Theorem~\ref{Thm:ApproxPropertie pp-pressure} below is preceded by the well-posedness of the discrete
problem~\eqref{discrStokes} for the pressure-modified Scott-Vogelius element that is a
consequence of the inf-sup stability inherited from the classical Scott-Vogelius element.
The following proposition recalls a right-inverse of the divergence from~\cite{Ainsworth_parker_hp_version} that is
bounded in terms of $\Theta_{\min}^{-1}$ with $\Theta_{\min}$ from~\eqref{thetamin}.
\begin{proposition}[{\cite{Ainsworth_parker_hp_version}}]
\label{Prop:InfSup stability mod SV} For $k\geq4$, there is a linear
operator $\Pi_{k}:M_{0,k-1}\left(  \mathcal{T}\right)  \rightarrow
\mathbf{S}_{k,0}\left(  \mathcal{T}\right)  $ with
\begin{equation}
\operatorname{div}\Pi_{k}q=q\quad\text{and\quad}\left\Vert \Pi_{k}q\right\Vert
_{\mathbf{H}^{1}\left(  \Omega\right)  }\leq\left(  c\Theta_{\min}\right)
^{-1}\left\Vert q\right\Vert _{L^{2}\left(  \Omega\right)  }\quad\forall q\in
M_{0,k-1}\left(  \mathcal{T}\right)  . \label{rightinverseest}%
\end{equation}
The constant $c$ only depends on $\Omega$ and the shape-regularity
of the mesh. The inf-sup constant for the Scott-Vogelius element $\left(
\mathbf{S}_{k,0}\left(  \mathcal{T}\right)  ,M_{0,k-1}\left(  \mathcal{T}%
\right)  \right)  $ is bounded from below by $c\Theta_{\min}$.\qed
\end{proposition}
Recall the inf-sup constant $\beta$ from~\eqref{infsupcond}, $C_f$ from Theorem~\ref{Thm:ApproxPropertie
pp-pressure}, and $c$ from Proposition~\ref{Prop:InfSup stability mod
SV}.
\begin{lemma}
\label{Propinfsupmod}Let $k\geq4$. The pressure-improved Scott-Vogelius
element 
is inf-sup stable with%
\begin{equation}
\label{Eq:SV infsup estimate}\beta\big(  \mathbf{S}_{k,0}\left(
\mathcal{T}\right)  ,M_{0,k-1}^{\operatorname{mod}}\left(  \mathcal{T}\right)
\big)  \geq c\Theta_{\min}/C_f.
\end{equation}
\end{lemma}

%

\begin{proof}
The operator $\mathcal{E}_{k}:M_{0,k-1}\left(  \mathcal{T} 
\right)  \rightarrow M_{0,k-1}^{\operatorname{mod}}\left(  \mathcal{T}\right)$ given for any 
$q\in M_{0,k-1}(\mathcal{T})$ by
\begin{equation}
\mathcal{E}_{k}q=q+\sum_{\mathbf{z}\in\mathcal{SC}_{\mathcal{T}}}%
f_{\mathbf{z}}\left(  q\right)  \left(  b_{k-1,\mathbf{z}}\right)
_{\operatorname*{mvz}} \label{defextraop}%
\end{equation}
is surjective onto $M_{0,k-1}^{\operatorname*{mod}}\left(  \mathcal{T}\right)$ by~\eqref{Eq:Def SV Mmod}.
Proposition \ref{PropOrthCompl}
reveals for any $q\in M_{0,k-1}(\mathcal{T})$ that
\begin{align}
\left(  \mathcal{E}_{k}q,q\right)
_{L^{2}\left(  \Omega\right)  }  &  =\|q\|_{L^{2}\left(  \Omega\right)  }^2+\sum_{\mathbf{z}%
\in\mathcal{SC}_{\mathcal{T}}}f_{\mathbf{z}}\left(  q\right)  \left(  \left(
b_{k-1,\mathbf{z}}\right)  _{\operatorname*{mvz}},q\right)  _{L^{2}\left(  \Omega\right)  }\label{expopest1}
=\|q\|_{L^{2}\left(  \Omega\right)  }^2.
\end{align}
Triangle inequalities, the boundedness of $f_{\mathbf{z}}$ from~\eqref{estfz}, and $C_f=1+\sum^{}_{\mathbf{z}\in
\mathcal{S}\mathcal{C}_{\mathcal{T}}} C_{f_{\mathbf{z}}}$ show
\[
\left\Vert \mathcal{E}_{k}q\right\Vert _{L^{2}\left(  \Omega\right)  }%
\leq\left\Vert q\right\Vert _{L^{2}\left(  \Omega\right)  }+\sum
_{\mathbf{z}\in\mathcal{SC}_{\mathcal{T}}}\left\vert f_{\mathbf{z}}\left(
q\right)  \right\vert \left\Vert b_{k-1,\mathbf{z}}\right\Vert _{L^{2}\left(
\omega_{\mathbf{z}}\right)  }\leq C_f  \left\Vert q\right\Vert
_{L^{2}\left(  \Omega\right)  }.%
\]
This, the surjectivity of $\mathcal{E}_k$, and the choice $\mathbf{v}%
_{q}:=\Pi_{k}q\in \mathbf{S}_{k,0}(\mathcal{T})$ with $\operatorname{div}(\mathbf{v}_q)=q$ and
$c\Theta_{\min}\|\mathbf{v}_q\|_{\mathbf{H}^1(\Omega)}\leq \|q\|_{L^2(\Omega)}$ from
Proposition~\ref{Prop:InfSup stability mod SV} for any $q\in M_{0,k-1}(\mathcal{T})$ verify
\begin{align*}
\inf_{q^{\ast}\in M_{0,k-1}^{\operatorname{mod}}\left(  \mathcal{T}\right)
\setminus\left\{  0\right\}  }\sup_{\mathbf{v}\in\mathbf{S}_{k,0}\left(
\mathcal{T}\right)  \setminus\left\{  \mathbf{0}\right\}  }\frac{\left(
q^{\ast},\operatorname{div}\mathbf{v}\right)  _{L^{2}\left(  \Omega\right)  }%
}{\left\Vert q^{\ast}\right\Vert _{L^{2}\left(  \Omega\right)  }\left\Vert
\mathbf{v}\right\Vert _{H^{1}\left(  \Omega\right)  }}  &  \geq\inf_{q\in
M_{0,k-1}\left(  \mathcal{T}\right)  \setminus\left\{  0\right\}  }%
\frac{\left(  \mathcal{E}_{k}q,q\right)
_{L^{2}\left(  \Omega\right)  }}{\left\Vert \mathcal{E}_{k}q\right\Vert
_{L^{2}\left(  \Omega\right)  }\left\Vert \mathbf{v}_{q}\right\Vert
_{H^{1}\left(  \Omega\right)  }}\\
&  \geq c\Theta_{\min}/C_f.\qedhere
\end{align*}%
\end{proof}

\begin{proof}[Proof of Theorem \ref{Thm:ApproxPropertie pp-pressure}]
	Given any 
	$\mathbf{v}\in\mathbf{S}_{k,0}\left(  \mathcal{T}\right)$, Proposition~\ref{PropOrthCompl} verifies as
	in~\eqref{expopest1} that
\[
b\left(  \mathbf{v},p^{\ast}_{M}\right)  =b\left(  \mathbf{v},p_M\right)
+\sum_{\mathbf{z}\in\mathcal{SC}_{\mathcal{T}}}f_{\mathbf{z}}\left(
p_M\right)  b\left(  \mathbf{v},\left(  b_{k-1,\mathbf{z}}\right)
_{\operatorname*{mvz}}\right)  =b\left(  \mathbf{v},p_M\right).
\]
Since $(\mathbf{u}_{\mathbf{S}},p_M)$ solves~\eqref{discrStokes} for $\mathbf{S}=\mathbf{S}_{k,0}(\mathcal{T})$ and
$M=M_{0,k-1}(\mathcal{T})$ with 
$\operatorname{div} \mathbf{u}_S = 0$, 
this verifies that $\left( \mathbf{u}_{\mathbf{S}}, p_M^{\ast} \right)$ solves~\eqref{discrStokes} for $\mathbf{S}=\mathbf{S}_{k,0}(\mathcal{T})$ and
$M=M_{0,k-1}^{\operatorname{mod}}(\mathcal{T})$.
The inf-sup stability from Proposition~\ref{Prop:InfSup stability mod SV} and Lemma~\ref{Propinfsupmod} verify the
uniqueness of
the respective discrete solutions.
The pressure estimate follows from \cite[Chap. 5,
Thm. 5.2.3, (5.2.27)]{BoffiBrezziFortin} in combination with \eqref{Eq:SV infsup estimate}.
\end{proof}

The $hp$-explicit estimates of the first infimum $\inf_{\mathbf{v}%
\in\mathbf{S}_{k,0}\left(  \mathcal{T}\right)  }\left\Vert \mathbf{u}%
-\mathbf{v}\right\Vert _{H^{1}\left(  \Omega\right)  }$ in
(\ref{Eq:QuasiBestApprox of Mmod}) for functions $\mathbf{u}$ with certain
Sobolev smoothness are well known from the literature on $hp$ finite elements.
Here, our focus will be on the estimate of the second infimum in \eqref{Eq:QuasiBestApprox of Mmod} related to our
new pressure space.

\subsection{Assumptions for and examples of functionals $f_{\mathbf{z}}$}\label{sub:assumptions_on_fz}

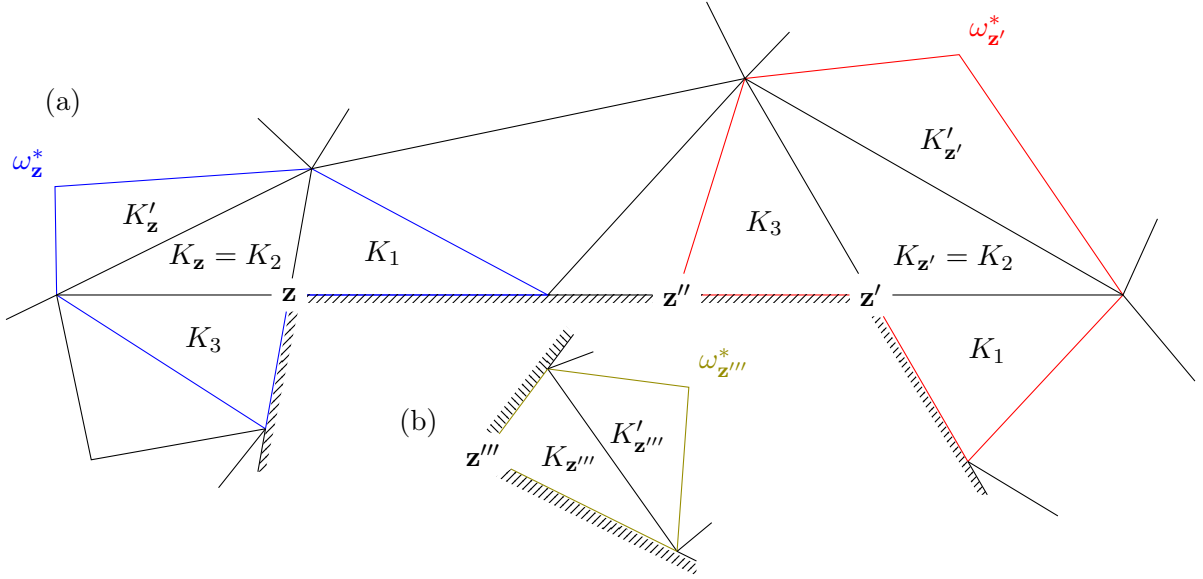
\begin{figure}[ptb]
\centering
\begin{tikzpicture}[scale=0.85]
\draw (0,0) node[] (z) {$\mathbf{z}$};
\coordinate (A) at (4,0);
\coordinate (A1) at (4.5,0);
\coordinate (B) at ($(z)!0.5!80:(A)$);
\coordinate (B1) at ($(z)!1.0!155:(A)$);
\coordinate (B2) at ($(z)!1.2!35:(A)$);
\coordinate (B3) at ($(z)!1.3!25:(A)$);
\coordinate (B4) at ($(z)!0.7!100:(A)$);
\coordinate (C) at ($(z)!0.9!180:(A)$);
\coordinate (C1) at ($(z)!1.0!220:(A)$);
\coordinate (C2) at ($(z)!1.1!185:(A)$);
\coordinate (D) at ($(z)!0.53!260:(A)$);
\coordinate (D1) at ($(z)!0.7!260:(A)$);
\coordinate (D2) at ($(z)!0.8!250:(A)$);
\draw (9,0) node[] (zprime) {$\mathbf{z}^{\prime}$};
\draw (6,0) node[] (E) {$\mathbf{z}^{\prime \prime}$};
\coordinate (F) at ($(zprime)!1!120:(E)$);
\coordinate (F1) at ($(zprime)!1.2!120:(E)$);
\coordinate (F2) at ($(zprime)!1.5!130:(E)$);
\coordinate (G) at ($(zprime)!1.3!180:(E)$);
\coordinate (G1) at ($(zprime)!1.33!250:(E)$);
\coordinate (G2) at ($(zprime)!1.53!195:(E)$);
\coordinate (G3) at ($(zprime)!1.73!165:(E)$);
\coordinate (H) at ($(zprime)!1.3!300:(E)$);
\coordinate (H1) at ($(H)!0.3!40:(G1)$);
\coordinate (H2) at ($(H)!0.4!110:(G1)$);
\coordinate (B5) at ($(B2)!0.4!50:(H)$);
\coordinate (B6) at ($(B2)!0.6!150:(H)$);
\coordinate (B7) at ($(B)!0.3!46:(B2)$);
\draw[blue] (z) -- (A) -- (B) -- (B1) -- (C) -- (D);
\draw (B) -- (C) -- (C1) -- (D);
\draw (A) -- (H);
\draw (B) -- (B4);
\draw[blue] (D) -- (z) -- (A);
\draw (z) -- (B);
\draw (z) -- (C) -- (C2);
\draw (D) -- (D2);
\draw (D) -- (D1);
\draw (B) -- (B7);
\draw (B)  --(H);
\draw (A) -- (E);
\draw[red] (E) -- (H) -- (G1) -- (G) --(F);
\draw (zprime) -- (H) -- (G);
\draw[red] (F) -- (zprime) -- (E);
\draw (F) -- (F1);
\draw (zprime) -- (G);
\draw (F) -- (F2);
\draw (G) -- (G2);
\draw (G) -- (G3);
\draw (H) -- (H1);
\draw (H) -- (H2);
\node (T1) at (barycentric cs:A=1,B=1,z=1) {$K_{1}$};
\node (T1) at (barycentric cs:C=1,D=1,z=1) {$K_{3}$};
\node (T1) at (barycentric cs:z=1.3,C=1,B=1) {$K_{\mathbf{z}} = K_2$};
\node (T1) at (barycentric cs:C=1,B=1,B1=1) {$K_{\mathbf{z}}^{\prime}$};
\node[above left] at (B1) {\textcolor{blue}{$\omega_{\mathbf{z}}^{\ast}$}};
\node[above right] at (G1) {\textcolor{red}{$\omega_{\mathbf{z}^{\prime}}^{\ast}$}};
\node (T1) at (barycentric cs:F=1,G=1,zprime=1) {$K_{1}$};
\node (T1) at (barycentric cs:E=1,H=1,zprime=1) {$K_{3}$};
\node (T1) at (barycentric cs:zprime=1.2,G=1.2,H=0.5) {$K_{\mathbf{z}^{\prime}} = K_2$};
\node (T1) at (barycentric cs:H=1,G=1,G1=1) {$K_{\mathbf{z}^{\prime}}^{\prime}$};
\fill[pattern=north east lines, ] ([xshift=0.00em]z.south) -- (D1) -- ([xshift=.3em]D1) -- ([xshift=.35em]z.south);
\fill[pattern=north east lines, ] (z.east) rectangle ([yshift=-.3em]E.west);
\fill[pattern=north east lines, ] (E.east) rectangle ([yshift=-.3em]zprime.west);
\fill[pattern=north east lines, ] ([xshift=0.35em]zprime.south) -- (F1) -- ([xshift=-.35em]F1) -- ([xshift=-.00em]zprime.south);

\draw (3,-2.5) node[] (z1) {$\mathbf{z}^{\prime \prime \prime}$};
\coordinate (a) at (6,-4);
\coordinate (a1) at ($(z1)!1.1!0:(a)$);
\coordinate (a2) at ($(z1)!1.1!10:(a)$);
\coordinate (b) at ($(z1)!0.5!80:(a)$);
\coordinate (b1) at ($(z1)!1!45:(a)$);
\coordinate (b2) at ($(z1)!0.7!80:(a)$);
\coordinate (b3) at ($(z1)!0.7!70:(a)$);
\coordinate (c) at ($(z1)!0.9!180:(a)$);
\coordinate (d) at ($(z1)!0.53!280:(a)$);
\draw (a) -- (b) -- (b2);
\draw[olive] (a) -- (b1) -- (b);
\draw[olive] (z1) -- (a);
\draw[olive] (z1) -- (b);
\draw (a) -- (a1); 
\draw (a) -- (a2);
\draw (b) -- (b3);
\phantom{		\node (T1) at (barycentric cs:z1=1,a=1,d=1) {$K_4$};
\node (T1) at (barycentric cs:z1=1,c=1,d=1) {$K_3$};}
\node (T1) at (barycentric cs:z1=1,a=1,b=1) {$K_{\mathbf{z}^{\prime \prime \prime}}$};
\node (T1) at (barycentric cs:b1=1,a=1,b=1) {$K_{\mathbf{z}^{\prime \prime \prime}}^{\prime}$};
\fill[pattern=north east lines, ] (z1) -- ([yshift=-.13em]a1) --([yshift=-.6em]a1) --([xshift=.52em]z1.south);
\fill[pattern=north west lines, ] ([xshift=0.58em]z1.north) -- (b2.north) -- ([xshift=-.3em,yshift=.32em]b2.north) -- ([xshift=.8,yshift=.10em]z1.north);
\node[above right] at (b1) {\textcolor{olive}{$\omega_{\mathbf{z}^{\prime \prime \prime}}^{\ast}$}};
\draw (-3.5,3) node[] (t) {(a)};
\draw (2,-2) node[] (t) {(b)};
\end{tikzpicture}
\caption{\justifying (a) The extended vertex patch $\omega_{\mathbf{z}}^{\ast}$ (contoured by blue lines) is of Robinson type, whereas the extended vertex patch $\omega_{\mathbf{z}^{\prime}}^{\ast}$ (outlined in red) is not of Robinson type due to presence of $\mathbf{z}^{\prime \prime} \in \mathcal{C}_{\mathcal{T}} \setminus \mathcal{SC}_{\mathcal{T}}$. Observe that $\mathbf{z}, \mathbf{z}^{\prime} \in \mathcal{SC}_{\mathcal{T}}$ satisfy \eqref{Item:seperated} in Def. \ref{Def:Robinson vertex}. (b) The extended vertex patch $\omega_{\mathbf{z}^{\prime \prime \prime}}^{\ast}$ is of Robinson type for $\mathbf{z}^{\prime \prime \prime} \in \mathcal{SC}_{\mathcal{T}}$ with $N_{\mathbf{z}^{\prime \prime \prime}} = 1$ (contoured in olive).} 
\label{Fig:Extendend SC nodal patches}%
\end{figure}

In order to derive approximation properties for the modified discrete pressure
space $M_{0,k-1}^{\operatorname{mod}}\left(  \mathcal{T}\right)  $, we will
specify the functionals $f_{\mathbf{z}}$ in a concrete way. In order to reduce
technicalities, we impose some simplifying assumptions on the
triangulation (see, Figure \ref{Fig:Extendend SC nodal patches}). 
Take any $\mathbf{z}\in\mathcal{SC}_{\mathcal{T}}$ with nodal patch
$\mathcal{T}_{\mathbf{z}}:=\left\{  K_{j}:1\leq j\leq N_{\mathbf{z}}\right\}
$ and set $K_{\mathbf{z}}\in\mathcal{T}_{\mathbf{z}}$ to be
\begin{equation}
K_{\mathbf{z}}:=%
\begin{cases}
K_{1} & \text{if }N_{\mathbf{z}}=1,\\
K_{2} & \text{if }N_{\mathbf{z}}=3.
\end{cases}
\label{Kzintro}%
\end{equation}

\begin{assumption}
\label{AExt}For any $\mathbf{z}\in\mathcal{SC}_{\mathcal{T}}$ with
$K_{\mathbf{z}}$ as in (\ref{Kzintro}), there exists a triangle denoted by
$K_{\mathbf{z}}^{\prime}\in\mathcal{T}$ which is adjacent to $K_{\mathbf{z}}$
but not contained in $\mathcal{T}_{\mathbf{z}}$; see Figure
\ref{Fig:Extendend SC nodal patches} for reference.
\end{assumption}

Assumption \ref{AExt} allows us to define the extended nodal patch
for $\mathbf{z}\in\mathcal{SC}_{\mathcal{T}}$ by
\[
\mathcal{T}_{\mathbf{z}}^{\ast}:=\mathcal{T}_{\mathbf{z}}\cup\left\{
K_{\mathbf{z}}^{\prime}\right\}  \quad\text{and\quad}\omega_{\mathbf{z}}%
^{\ast}:=\omega_{\mathbf{z}}\cup K_{\mathbf{z}}^{\prime}.
\]
To reduce technicalities, we restrict to meshes where super-critical vertices
are properly separated. The definition is illustrated by Fig.
\ref{Fig:Extendend SC nodal patches}.

\begin{definition}
\label{Def:Robinson vertex} Let Assumption \ref{AExt} hold. A super-critical vertex $\mathbf{z}\in
\mathcal{SC}_{\mathcal{T}}$ is called a \emph{Robinson vertex }if it is
\emph{isolated}:

\begin{enumerate}
\item \label{Item:seperated}$\mathcal{T}_{\mathbf{z}}^{\ast}\cap\mathcal{T}_{\mathbf{y}}^{\ast
}=\emptyset$ for all $\mathbf{y}\in\mathcal{SC}_{\mathcal{T}}\setminus\left\{
\mathbf{z}\right\}  $.

\item $\omega_{\mathbf{z}}^{\ast}\cap\mathcal{C}_{\mathcal{T}}=\left\{
\mathbf{z}\right\}  $.
\end{enumerate}
\end{definition}

\begin{remark}
\label{Rem:Avoiding techinacl difficulties}
We assume below that all super critical vertices are of Robinson type to reduce technicalities:
If $\mathbf{z}\in\mathcal{S}\mathcal{C}_{\mathcal{T}}$ is a Robinson vertex,
the critical function $b_{k-1,\mathbf{z}}$ satisfies
\begin{enumerate}
	\item[(i)] $(b_{k-1,\mathbf{z}}, b_{k-1,\mathbf{y}})_{L^2(\Omega)}=0$ for all
		$\mathbf{y}\in\mathcal{C}_{\mathcal{T}}\setminus\{\mathbf{z}\}$ and
	\item[(ii)] $b_{k-1,\mathbf{z}}|_{K_\mathbf{y}'}=0$ for all $\mathbf{y}\in \mathcal{S}\mathcal{C}_{\mathcal{T}}$.
\end{enumerate}
These properties 
allow us to investigate the approximation property locally in a vicinity of
Robinson vertices and avoid clusters of super critical vertices with additional
coupling effects.
\end{remark}

The following \emph{structural} assumption on 
$f_{\mathbf{z}}:\mathbb P_{k-1}(\mathcal{T})\to \mathbb R$ from~\eqref{estfz}
leaves considerable freedom for the particular choice and enables optimal approximation properties of
$M_{0,k-1}^{\operatorname{mod}}(\mathcal{T})$ in Theorem~\ref{Thm:Recoverd bestapprox} below.
For a polynomial $q\in\mathbb{P}_{k}\left(  K\right)  $ we
denote its analytic extension to $\mathbb{R}^{2}$ by $q^{\operatorname*{ext}}$. 


\begin{assumption}
\label{Astruct}
Let with Assumption \ref{AExt} hold and
$f_{\mathbf{z}}:\mathbb P_{k-1}(\mathcal{T})\to\mathbb R$ be given in terms of continuous linear functionals $J_{\mathbf{z}}:\mathbb P_{k-1}(\mathbb
R^2)\to\mathbb R$ by
\begin{equation}
f_{\mathbf{z}}\left(  q\right)  :=J_{\mathbf{z}}\left(  \left.  q\right\vert
_{K_{\mathbf{z}}^{\prime}}^{\operatorname{ext}}\right)  -J_{\mathbf{z}}\left(
\left.  q\right\vert _{K_{\mathbf{z}}}^{\operatorname{ext}}\right) 
\label{formfz}%
\end{equation}
for any $\mathbf{z}\in\mathcal{S}\mathcal{C}_{\mathcal{T}}$. 
Each of the functionals $J_\mathbf{z}$ satisfies
\begin{equation}
J_{\mathbf{z}}\left(  \left.  b_{k-1,\mathbf{z}}\right\vert _{K_{\mathbf{z}}%
}^{\operatorname{ext}}\right)  =1\quad\text{and}\quad
\left\vert J_{\mathbf{z}}\left(  q\right)  \right\vert \leq\frac
{C_{\mathbf{z}}}{\left\Vert b_{k-1,\mathbf{z}}\right\Vert _{L^{2}\left(
\omega_{\mathbf{z}}\right)  }}\left\Vert q\right\Vert _{L^{2}\left(
\mathcal{U}_{\mathbf{z}}\right)  }\quad\forall q\in\mathbb{P}_{k-1}\left(
\mathbb{R}^{2}\right)  \label{Eq:Jzcontestimate}%
\end{equation}
for some constant $C_{\mathbf{z}}$ which is independent of $q$ and $h_{\mathbf{z}}$ but, possibly, depends on $k$, $\Omega$, and $\gamma_{\mathcal{T}}$ and some subset %
$\mathcal{U}_{\mathbf{z}}\subset K_{\mathbf{z}}\cup K_{\mathbf{z}}^{\prime}$ 
with the following property.
There exists some $\delta_{\mathbf{z}}\geq0$ and triangles $K_{\mathbf{z}}%
^{\operatorname*{ext}}$, $K_{\mathbf{z}}^{\prime,\operatorname*{ext}}%
\subseteq\mathbb{R}^{2}$ not necessarily contained in $\mathcal{T}$ or
$\Omega$ such that%
\begin{subequations}
\label{Kzext}
\end{subequations}%
\begin{align}
K_{\mathbf{z}}\cup\mathcal{U}_{\mathbf{z}}  &  \subset K_{\mathbf{z}%
}^{\operatorname*{ext}}\subset\left\{  \mathbf{y}\in\mathbb{R}^{2}%
\mid\operatorname{dist}\left(  \mathbf{y},K_{\mathbf{z}}\right)  \leq
\delta_{\mathbf{z}}h_{K_{\mathbf{z}}}\right\}  ,\tag{%
\ref{Kzext}%
a}\label{Kzexta}\\
K_{\mathbf{z}}^{\prime}\cup\mathcal{U}_{\mathbf{z}}  &  \subset K_{\mathbf{z}%
}^{\prime,\operatorname*{ext}}\subset\left\{  \mathbf{y}\in\mathbb{R}^{2}%
\mid\operatorname{dist}\left(  \mathbf{y},K_{\mathbf{z}}^{\prime}\right)
\leq\delta_{\mathbf{z}}h_{K_{\mathbf{z}}^{\prime}}\right\}.  \tag{%
\ref{Kzext}%
b}\label{Kzextb}%
\end{align}

\end{assumption}

Next we introduce an overlap constant for some local neighborhoods of super-critical
vertices. For $\mathbf{z}\in\mathcal{SC}_{\mathcal{T}}$ define the regions
(cf. (\ref{defomegaK}))%
\begin{equation}
\widetilde{\omega}_{\mathbf{z}}:=\omega\left(  K_{\mathbf{z}}\right)
\cup\omega\left(  K_{\mathbf{z}}^{\prime}\right)  \cup K_{\mathbf{z}%
}^{\operatorname*{ext}}\cup K_{\mathbf{z}}^{\prime,\operatorname*{ext}}.
\label{defomegazstar}%
\end{equation}
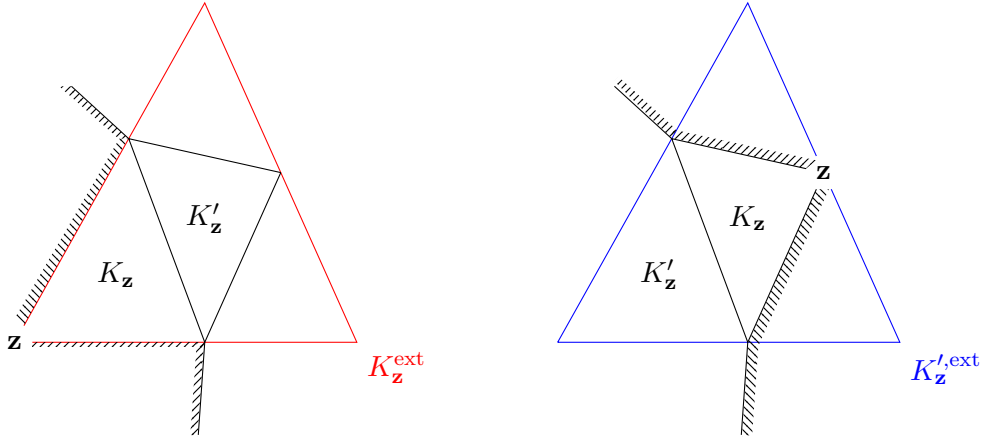
\begin{figure}[ptb]
\centering\begin{tikzpicture}[scale=0.5]
\coordinate (z) at (7,4.5);
\coordinate (B) at (9,0);
\coordinate (C) at (5,9);
\draw (0,0) node[] (A) {$\mathbf{z}$};
\coordinate (z1) at (3,27/5);
\coordinate (z2) at (5,0);
\coordinate (y1) at ($(z1)!0.5!150:(z)$);
\coordinate (y2) at ($(z2)!0.5!200:(z)$);

\draw[red] (A) -- (B) -- (z) -- (C) -- (A);
\draw (z) -- (z1) -- (z2) -- (z);
\draw (z1) -- (y1);
\draw (z2) -- (y2);
\fill[pattern=north east lines, ] (A.east) rectangle ([yshift=-.3em]z2.west);
\fill[pattern=north west lines, ] ([xshift=0.0em]A) -- ([xshift=-.05em]z1.west) -- ([xshift=-.65em]z1.west) -- ([yshift=.50em]A.north);
\fill[pattern=north east lines, ] ([xshift=-0.0em]z1) -- ([xshift=-.05em]y1.north) -- ([xshift=-.85em]y1.west) -- ([xshift=-.70em]z1.west);
\fill[pattern=north east lines, ] ([xshift=-.05em]z2.south) -- ([xshift=-.05em,yshift=-1]y2.south) -- ([xshift=-.55em]y2.south) -- ([xshift=-.55em]z2.south);

\node (T1) at (barycentric cs:A=1,z1=1,z2=1) {$K_{\mathbf{z}}$};
\node (T1) at (barycentric cs:z=1,z1=1,z2=1) {$K_{\mathbf{z}}^{\prime}$};
\node[below right] at (B) {\textcolor{red}{$K_{\mathbf{z}}^{ \operatorname{ext}}$}};
\end{tikzpicture}
\hfil
\begin{tikzpicture}[scale=0.5]
\coordinate (A) at (0,0);
\coordinate (B) at (9,0);
\coordinate (C) at (5,9);
\draw (7,4.5) node[] (z) {$\mathbf{z}$};
\coordinate (z1) at (3,27/5);
\coordinate (z2) at (5,0);
\coordinate (y1) at ($(z1)!0.5!150:(z)$);
\coordinate (y2) at ($(z2)!0.5!200:(z)$);

\draw[blue] (A) -- (B) -- (z) -- (C) -- (A);
\draw (z) -- (z1) -- (z2) -- (z);
\draw (z1) -- (y1);
\draw (z2) -- (y2);
\fill[pattern=north west lines, ] ([xshift=.05em]z) -- ([xshift=.05em]z2.north) -- ([xshift=.75em,yshift=-.2em]z2.east) -- ([xshift=.35em]z.south);
\fill[pattern=north east lines, ] ([xshift=0.0em]z) -- ([yshift=.05em]z1) -- ([yshift=.65em]z1) -- ([xshift=-.50em,yshift=-.30em]z.north);
\fill[pattern=north east lines, ] ([xshift=0.0em]z1.north) -- ([yshift=.35em]y1.north) -- ([xshift=.3em,yshift=.65em]y1.east) -- ([xshift=.3em,yshift=.650em]z1.north);
\fill[pattern=north west lines, ] ([xshift=.05em]z2.south) -- ([xshift=.05em,yshift=-1]y2.south) -- ([xshift=.75em]y2.south) -- ([xshift=.75em]z2.south);

\node (T1) at (barycentric cs:z=1,z1=1,z2=1) {$K_{\mathbf{z}}$};
\node (T1) at (barycentric cs:A=1,z1=1,z2=1) {$K_{\mathbf{z}}^{\prime}$};
\node[below right] at (B) {\textcolor{blue}{$K_{\mathbf{z}}^{\prime, \operatorname{ext}}$}};
\end{tikzpicture}

\caption{Ilustration of the extended triangles $K_{\mathbf{z}}^{\operatorname{ext}}$ (left in red) and $K_{\mathbf{z}}^{\prime ,\operatorname{ext}}$ (right in blue) for the case $\mathbf{z} \in \mathcal{SC}_{\mathcal{T}}$ with $N_{\mathbf{z}} = 1$.}%
\label{Fig:Extendend triangles}%
\end{figure}
The maximal overlap is described by
\begin{equation}
\max_{\mathbf{z}\in\mathcal{SC}_{\mathcal{T}}}\operatorname{card}\left\{
\mathbf{y}\in\mathcal{SC}_{\mathcal{T}}\mid\left\vert \widetilde{\omega
}_{\mathbf{z}}\cap\widetilde{\omega}_{\mathbf{y}}\right\vert >0\right\}
=:C_{\operatorname*{ov}}. \label{defCo}%
\end{equation}
Note that it is a very mild assumption on the mesh to assume that $C_{\operatorname{ov}}$ is moderately bounded; all constants in the sequel may depend on $C_{\operatorname*{ov}}$. Note that
the constants
\begin{equation}
C_{J}:=\max_{\mathbf{z}\in\mathcal{SC}_{\mathcal{T}}}C_{\mathbf{z}}%
\quad\text{and\quad}\delta_{\max}:=\max_{\mathbf{z}\in\mathcal{SC}%
_{\mathcal{T}}}\delta_{\mathbf{z}} \label{CJconstant}%
\end{equation}
are independent of $\mathbf{z}\in\mathcal{SC}_{\mathcal{T}}$ and $h_{\mathbf{z}}$ but,
possibly, depend on $k$, $\Omega$ and $\gamma_{\mathcal{T}}$ (see Assumption \ref{Astruct}).

\subsection{The approximation property for the improved pressure space
$M_{0,k-1}^{\operatorname*{mod}}\left(  \mathcal{T}\right)  $%
\label{ApproxSVpressure}}

This subsection verifies the optimal approximation properties of the
modified pressure space $M_{0,k-1}^{\operatorname{mod}} \left(  \mathcal{T}
\right)$, \eqref{Eq:Def SV Mmod} under
Assumption \ref{Astruct} with explicit constants in terms of the mesh width.
Consider the continuous Stein extension 
$\mathcal{E}_{\operatorname{Stein}}:C^{\infty}\left(  \overline{\Omega}\right)
\rightarrow C^{\infty}\left(  \mathbb{R}^{2}\right)$ from \cite[Thm. 5, p.~181]{emstein} that extends, for any $m\geq0$, to a
continuous operator $\mathcal{E}_{\operatorname{Stein}}:H^{m}\left(
\Omega\right)  \rightarrow H^{m}\left(  \mathbb{R}^{2}\right)$ with
\begin{equation}
\left\Vert \mathcal{E}_{\operatorname{Stein}}u\right\Vert _{H^{m}\left(
\mathbb{R}^{2}\right)  }\leq C_{\operatorname{Stein}}\left\Vert u\right\Vert
_{H^{m}\left(  \Omega\right)  }\quad\forall u\in H^{m}\left(  \Omega\right).
\label{Cstein}%
\end{equation}
The constant $C_{\operatorname{Stein}%
}$ depends on $m$ and $\Omega$. 
Let $T_{k}\in \mathbb P_k(\mathbb R)$ denote the Chebyshev polynomial of
first kind and degree $k$; see \cite[Table 18.3.1]{NIST:DLMF}
for details.
\begin{theorem}
\label{Thm:Recoverd bestapprox} Let Assumption \ref{Astruct} hold and suppose that
all super-critical vertices are of Robinson type. 
For any $p\in H^{s-1}\left(
\Omega\right)  \cap L_{0}^{2}\left(  \Omega\right)  $ with $s>1$,
there exists $p_{M}\in M_{0,k-1}^{\operatorname*{mod}}\left(  \mathcal{T}%
\right)  $ such that%
\begin{equation}
\left\Vert p-p_{M}\right\Vert _{L^{2}\left(  \mathbb{R}^{2}\right)  }\leq
C_{\operatorname*{apx}}\frac{\left(  \left(  1+2\delta_{\max}\right)
h_{\mathcal{T}}\right)  ^{\min\left\{  k,s-1\right\}  }}{k^{s-1}}\left\Vert
p\right\Vert _{H^{s-1}\left(  \Omega\right)  } \label{errorest}%
\end{equation}
holds with $C_{\operatorname*{apx}}:=C\sqrt{\operatorname{card}\mathcal{SC}_{\mathcal{T}}%
}C_{\operatorname{Stein}}C_{J}T_{k-1}\left(  1+c\delta_{\max}\right)$.
The constant $C$ depends only on $C_{\operatorname*{ov}}$ from~\eqref{defCo}
and on the shape regularity of the mesh.
\end{theorem}

An auxiliary result on the decomposition of functions 
$\widetilde{Q}_{0}^{h,k}$
from \eqref{defQtilde} precedes 
the proof of Theorem \ref{Thm:Recoverd bestapprox} below.

\begin{lemma}
\label{Lem:Decomposition of Q0hk}Let all super critical vertices $\mathbf{z}
\in\mathcal{SC}_{\mathcal{T}}$ be of Robinson type and let $k\geq4$ be given.
Then for all $\widetilde{q}\in\widetilde{Q}_{0}^{h,k}$ there exists $q\in
M_{0,k-1}\left(  \mathcal{T}\right)  $ and $\theta_{\mathbf{z}}\in\mathbb{R}$
for all $\mathbf{z}\in\mathcal{SC}_{\mathcal{T}}$ such that $\tilde{q}$ can be
written as
\[
\widetilde{q}=q+\sum_{\mathbf{z}\in\mathcal{SC}_{\mathcal{T}}}\theta
_{\mathbf{z}}\left(  b_{k-1,\mathbf{z}}\right)  _{\operatorname*{mvz}}.
\]

\end{lemma}%

\begin{proof}
Let $\widetilde{q}\in\widetilde{Q}_{0}^{h,k}\subset \mathbb P_{k-1,0}(\mathcal{T})$ be arbitrarily chosen but fixed.
The orthogonal decomposition of $\mathbb P_{k-1,0}(\mathcal{T})$ in Proposition \ref{PropOrthCompl} provides $q\in
M_{0,k-1}\left(  \mathcal{T}\right)  $ and $\theta_{\mathbf{z}}\in\mathbb R$ for all 
$\mathbf{z}\in\mathcal{C}_{\mathcal{T}}$ such that
\[
\widetilde{q}=q+\sum_{\mathbf{z}\in\mathcal{C}_{\mathcal{T}}}\theta
_{\mathbf{z}}\left(  b_{k-1,\mathbf{z}}\right)  _{\operatorname*{mvz}}.
\]
For $\mathbf{z}\in\mathcal{C}_{\mathcal{T}}\setminus\mathcal{SC}_{\mathcal{T}%
}$, we locally decompose $\mathbb{P}_{k-1}\left(  \mathcal{T}_{\mathbf{z}%
}\right)  $ into
\[
\mathbb{P}_{k-1}\left(  \mathcal{T}_{\mathbf{z}}\right)  =M_{0,k-1}\left(
\mathcal{T}_{\mathbf{z}}\right)  \oplus\operatorname{span}\left\{
b_{k-1,\mathbf{z}}\right\}  ,
\]
where $M_{0,k-1}\left(  \mathcal{T}_{\mathbf{z}}\right)  :=\left\{  \left.
q\in\mathbb{P}_{k-1}\left(  \mathcal{T}_{\mathbf{z}}\right)  \;\right\vert
\;A_{\mathcal{T},\mathbf{z}}\left(  q\right)  =0\right\}  $. Since $\tilde{q}$
is continuous at $\mathbf{z}$ and $N_{\mathbf{z}}$ is even since
$\mathbf{z}\in\mathcal{C}_{\mathcal{T}}\setminus\mathcal{SC}_{\mathcal{T}}$,
the term $A_{\mathcal{T},\mathbf{z}}\left(  q\right)  =0$ vanishes and
$\left.  q\right\vert _{\omega_{\mathbf{z}}}\in M_{0,k-1}\left(
\mathcal{T}_{\mathbf{z}}\right)  $ follows. This and $A_{\mathcal{T}%
,\mathbf{z}}\left(  \left(  b_{k-1,\mathbf{z}}\right)  _{\operatorname*{mvz}%
}\right)  \neq0$ from \eqref{Azpropb} implies $\theta_{\mathbf{z}}=0$
concluding the proof.
\end{proof}

This lemma enables the proof of the approximation property for the modified pressure space.

\begin{proof}[Proof of Theorem \ref{Thm:Recoverd bestapprox}]
Let $p\in H^{s-1}\left(  \Omega\right)  \cap L_{0}^{2}\left(  \Omega\right)  $ be given. With a
slight abuse of notation its extension by $\mathcal{E}_{\operatorname{Stein}}$
to $\mathbb{R}^{2}$ is again denoted by $p$.
Theorem 2.1 from \cite{Ainsworth_parker_II} provides some $p_{Q}\in\widetilde{Q}_{0}^{h,k}\left(  \mathcal{T}\right)  $
that satisfies, for any $K\in\mathcal{T}$, the estimate
\begin{equation} \label{Eq:Cont best aprox}
\left\Vert p-p_{Q}\right\Vert _{L^{2}\left(  K\right)  }\leq C\frac
{h_{K}^{\min\left\{  k,s-1\right\}  }}{k^{s-1}}\left\Vert p\right\Vert
_{H^{s-1}\left(  \omega\left(  K\right)  \right)  }.
\end{equation}
Lemma \ref{Lem:Decomposition of Q0hk} reveals, for some $q\in M_{0,k-1}\left(  \mathcal{T}\right)  $, $C_0 \in \mathbb{R}$ and $\theta_{\mathbf{z}%
}\in\mathbb{R}$ for all $\mathbf{z}\in\mathcal{SC}_{\mathcal{T}}$, the form
\begin{equation} \label{Eq:Decomposition cont best approx}
p_{Q}=q  +\sum_{\mathbf{z}\in\mathcal{SC}_{\mathcal{T}}}\theta_{\mathbf{z}%
}\left(  b_{k-1,\mathbf{z}}\right)  _{\operatorname*{mvz}}
=
q + C_{0}+ 
\sum_{\mathbf{z}\in\mathcal{SC}_{\mathcal{T}}}\theta_{\mathbf{z}}b_{k-1,\mathbf{z}%
}.
\end{equation}
Triangle inequalities result with $p_{M}:=q+\sum_{\mathbf{z}\in\mathcal{SC}_{\mathcal{T}}}f_{\mathbf{z}}\left(
q\right)  \left(  b_{k-1,\mathbf{z}}\right)_{\operatorname*{mvz}}\in M_{0,k-1}^{\operatorname{mod}} \left( \mathcal{T} \right)$ in 
\begin{equation}
\left\Vert p_{Q}-p_{M}\right\Vert _{L^{2}\left(  \Omega\right)  }\leq
\sum_{\mathbf{z}\in\mathcal{SC}_{\mathcal{T}}}\left\vert \theta_{\mathbf{z}%
}-f_{\mathbf{z}}\left(  q\right)  \right\vert \left\Vert b_{k-1,\mathbf{z}%
}\right\Vert _{L^{2}\left(  \Omega\right)  }. \label{phmpm}%
\end{equation}
Since the super-critical vertices $\mathbf{z} \in \mathcal{SC}_{\mathcal{T}}$ are of Robinson type
and Assumption \ref{AExt} holds, \eqref{Eq:Decomposition cont best approx} and Remark~\ref{Rem:Avoiding techinacl
difficulties}(ii) reveal $q\vert _{K_{\mathbf{z}}^{\prime}}=p_{Q}\vert
_{K_{\mathbf{z}}^{\prime}}-C_{0}$.
Hence,~\eqref{formfz} and algebraic manipulations result~in
\begin{align*}
\left\vert \theta_{\mathbf{z}}-f_{\mathbf{z}}\left(  q_{0}\right)
\right\vert  &  =\left\vert \theta_{\mathbf{z}}-\left(  J_{\mathbf{z}}\left(
\left.  p_{Q}\right\vert _{K_{\mathbf{z}}^{\prime}}^{\operatorname*{ext}%
}-C_{0}\right)  -J_{\mathbf{z}}\left(  \left.  \left(  p_{Q}-\theta
_{\mathbf{z}}b_{k-1,\mathbf{z}}\right)  \right\vert _{K_{\mathbf{z}}%
}^{\operatorname*{ext}}-C_{0}\right)  \right)  \right\vert \\
&  =\left\vert \theta_{\mathbf{z}}-J_{\mathbf{z}}\left(  \left.
p_{Q}\right\vert _{K_{\mathbf{z}}^{\prime}}^{\operatorname*{ext}}\right)
+J_{\mathbf{z}}\left(  \left.  p_{Q}\right\vert _{K_{\mathbf{z}}%
}^{\operatorname*{ext}}\right)  -\theta_{\mathbf{z}}J_{\mathbf{z}}\left(
\left.  b_{k-1,\mathbf{z}}\right\vert _{K_{\mathbf{z}}}^{\operatorname*{ext}%
}\right)  \right\vert
\overset{\text{(\ref{Eq:Jzcontestimate})}}{=}\left\vert J_{\mathbf{z}}\left(  \left.
p_{Q}\right\vert _{K_{\mathbf{z}}}^{\operatorname*{ext}}-\left.
p_{Q}\right\vert _{K_{\mathbf{z}}^{\prime}}^{\operatorname*{ext}}\right)
\right\vert .
\end{align*}
This, the bound \eqref{Eq:Jzcontestimate} of the functional $J_{\mathbf{z}}$ 
in Assumption~\ref{Astruct}, and triangle inequalities imply
\begin{align*}
C_{\mathbf{z}}^{-1}\left\Vert b_{k-1,\mathbf{z}}\right\Vert _{L^{2}\left(
\omega_{\mathbf{z}}\right)}&\left\vert \theta_{\mathbf{z}}-f_{\mathbf{z}}\left(  q_{0}\right)
\right\vert 
\leq \left\Vert p_{Q}-\left.
p_{Q}\right\vert _{K_{\mathbf{z}}^{\prime}}^{\operatorname*{ext}}\right\Vert
_{L^{2}\left(  \mathcal{U}_{\mathbf{z}}\right)  }+\left\Vert \left.
p_{Q}\right\vert _{K_{\mathbf{z}}}^{\operatorname*{ext}}-p_{Q}\right\Vert
_{L^{2}\left(  \mathcal{U}_{\mathbf{z}}\right)  }\\
&  \leq 2\left\Vert p_{Q}%
-p\right\Vert _{L^{2}\left(  K_{\mathbf{z}}\cup K_{\mathbf{z}}^{\prime
}\right)  }+\left\Vert \left.  p_{Q}\right\vert _{K_{\mathbf{z}}%
}^{\operatorname*{ext}}-p\right\Vert _{L^{2}\left(  K_{\mathbf{z}%
}^{\operatorname*{ext}}\right)  }+\left\Vert \left.  p_{Q}\right\vert
_{K_{\mathbf{z}}^{\prime}}^{\operatorname*{ext}}-p\right\Vert _{L^{2}\left(
K_{\mathbf{z}}^{\prime,\operatorname*{ext}}\right)  }
\end{align*}
with 
$\mathcal{U}_{\mathbf{z}} \subseteq (K_{\mathbf{z}} \cup K_{\mathbf{z}}^{\prime})\cap
K_{\mathbf{z}}^{\operatorname*{ext}}\cap K_{\mathbf{z}}^{\prime,\operatorname*{ext}}$ in the last step.
Lemma \ref{LemNeighborhood} with $\kappa\coloneqq\min\{k,s-1\}$ and the definition of the neighborhood $\widetilde{\omega}_{\mathbf{z}}$ in
\eqref{defomegazstar} control the second summand by
\begin{align*}
\left\Vert \left.  p_{Q}\right\vert _{K_{\mathbf{z}}}^{\operatorname*{ext}%
}-p\right\Vert _{L^{2}\left(  K_{\mathbf{z}}^{\operatorname*{ext}}\right)  }
&  \leq T_{k-1}\left(  1+c\delta_{\max}\right)  \left(2\frac{\left(  \left(
1+2\delta_{\max}\right)
h_{\mathbf{z}}\right)^{\kappa}}{k^{s-1}}\left\Vert p\right\Vert _{H^{s-1}\left(  K_{\mathbf{z}%
}^{\operatorname*{ext}}\right)  }
  +\left\Vert p-p_{Q}\right\Vert
_{L^{2}\left(  K_{\mathbf{z}}\right)  }\right)\\
&  \leq CT_{k-1}\left(  1+c\delta_{\max}\right)  \frac{\left(  \left(
1+2\delta_{\max}\right)  h_{\mathbf{z}}\right)  ^{\kappa
}}{k^{s-1}}\left\Vert p\right\Vert _{H^{s-1}\left(  \widetilde{\omega}_{\mathbf{z}} \right)  },
\end{align*}
with the approximation property (\ref{Eq:Cont best aprox}) of $p_Q$ in the last step.
The previous two estimates, the analogous estimate for $\left\Vert \left.  p_{Q}\right\vert _{K_{\mathbf{z}%
}^{\prime}}^{\operatorname*{ext}}-p\right\Vert _{L^{2}\left(  K_{\mathbf{z}%
}^{\prime, \operatorname{ext}}\right)  }$ with the same upper bound, and \eqref{Eq:Cont best aprox} result~in
\[
\left\vert \theta_{\mathbf{z}}-f_{\mathbf{z}}\left(  q_{0}\right)
\right\vert 
\leq CC_{\mathbf{z}}T_{k-1}\left(  1+c\delta_{\max}\right)  \frac{\left(  \left(
1+2\delta_{\max}\right)  h_{\mathbf{z}}\right)  ^{\kappa}}
{k^{s-1}\left\Vert
b_{k-1,\mathbf{z}}\right\Vert _{L^{2}\left(  \omega_{\mathbf{z}}\right)
}}\left\Vert p\right\Vert _{H^{s-1}\left(
\widetilde{\omega}_{\mathbf{z}}\right)  }.
\]
Hence, \eqref{phmpm} and $C_{\mathbf{z}}\leq C_J$ by definition in~\eqref{CJconstant} verify
\begin{align*}
\left\Vert p_{Q}-p_{M}\right\Vert _{L^{2}\left(  \Omega\right)  }  &  \leq
CC_{J}T_{k-1}\left(  1+c\delta_{\max}\right)  \frac{\left(  \left(
1+2\delta_{\max}\right)  h_{\mathcal{T}}\right)  ^{\kappa  }}{k^{s-1}} \sum_{\mathbf{z}\in\mathcal{SC}_{\mathcal{T}}}\left\Vert p\right\Vert
_{H^{s-1}\left( \widetilde{\omega}_{\mathbf{z}}\right)  }.
\end{align*}
A Cauchy inequality in $\ell^2$ and the finite overlay \eqref{defCo} of the regions $\widetilde{\omega}_{\mathbf{z}}$ implies%
\begin{align*}
\left\Vert p_{Q}-p_{M}\right\Vert _{L^{2}\left(  \Omega\right)  }  &  \leq
CC_{J}T_{k-1}\left(  1+c\delta_{\max}\right)  \sqrt{\operatorname{card}%
\mathcal{SC}_{\mathcal{T}}}\frac{\left(  \left(  1+2\delta_{\max}\right)
h_{\mathcal{T}}\right)  ^{\kappa}}{k^{s-1}}\left\Vert
p\right\Vert _{H^{s-1}\left(  \mathbb{R}^{2}\right)  }.
\end{align*}
This and $\|p\|_{H^{s-1}(\mathbb R^2)}\leq C_{\operatorname*{Stein}}\|p\|_{H^{s-1}(\Omega)}$ from the Stein extension (cf., \eqref{Cstein} conclude the proof.
\end{proof}

\subsection{Recovery of optimal rates for the postprocessed pressure}

\label{ApproxSVpressureII}

The structural assumptions on the functional $f_{\mathbf{z}}$ in Assumption~\ref{Astruct} for
optimal approximation properties of the resulting modified pressure space $M_{0,k-1}^{\operatorname{mod}%
}\left(  \mathcal{T}\right)  $ in Theorem~\ref{Thm:Recoverd bestapprox} leave considerable freedom in the particular
choice of the functionals
$f_{\mathbf{z}}$. Below we give two examples.

\begin{example}
\label{ExJ}Two possible choices of $J_{\mathbf{z}}$ and $f_{\mathbf{z}}$ are
given by

\begin{enumerate}
	\item the point evaluation at $\mathbf{z}\in \mathcal{S}\mathcal{C}_{\mathcal{T}}$, namely
\begin{equation}
J_{\mathbf{z}}\left(  q\right)  :=\frac{q\left(  \mathbf{z}\right)  }{\left.
b_{k-1,\mathbf{z}}\right\vert _{K_{\mathbf{z}}}\left(  \mathbf{z}\right)
}\quad\forall q\in\mathbb{P}_{k-1}\left(  \mathbb{R}^{2}\right),
\label{Jz0primi}%
\end{equation}
so that%
\begin{equation}
f_{\mathbf{z}}\left(  q\right)  =\frac{\left(  \left.  q\right\vert
_{K_{\mathbf{z}}^{\prime}}^{\operatorname*{ext}}-\left.  q\right\vert
_{K_{\mathbf{z}}}\right)  \left(  \mathbf{z}\right)  }{\left.
b_{k-1,\mathbf{z}}\right\vert _{K_{\mathbf{z}}}\left(  \mathbf{z}\right)
}\quad\forall q\in M_{0,k-1}\left(  \mathcal{T}\right)  , \label{Jzprimi}%
\end{equation}

\item the integration with weight 
	$b_{k-1,\mathbf{z}}\vert _{K_{\mathbf{z}}}^{\operatorname*{ext}}$ over 
	a subset $S_{\mathbf{z}}^{\prime}\subset K_{\mathbf{z}}^{\prime}$
	with positive measure, namely %
\[
J_{\mathbf{z}}\left(  q\right)  :=\frac{\int_{S_{\mathbf{z}}^{\prime}}q\left.
b_{k-1,\mathbf{z}}\right\vert _{K_{\mathbf{z}}}^{\operatorname*{ext}}%
}{\left\Vert \left.  b_{k-1,\mathbf{z}}\right\vert _{K_{\mathbf{z}}%
}^{\operatorname*{ext}}\right\Vert _{L^{2}\left(  S_{\mathbf{z}}^{\prime
}\right)  }^{2}}\quad\forall q\in\mathbb{P}_{k-1}\left(  \mathbb{R}%
^{2}\right),
\]
so that%
\begin{equation}
f_{\mathbf{z}}\left(  q\right)  =\frac{\int_{S_{\mathbf{z}}^{\prime}}\left(
q-\left.  q\right\vert _{K_{\mathbf{z}}}^{\operatorname*{ext}}\right)  \left.
b_{k-1,\mathbf{z}}\right\vert _{K_{\mathbf{z}}}^{\operatorname*{ext}}%
}{\left\Vert \left.  b_{k-1,\mathbf{z}}\right\vert _{K_{\mathbf{z}}%
}^{\operatorname*{ext}}\right\Vert _{L^{2}\left(  S_{\mathbf{z}}^{\prime
}\right)  }^{2}}\quad\forall q\in M_{0,k-1}\left(  \mathcal{T}\right)  .
\label{Szprime}%
\end{equation}
\end{enumerate}
The implementation of the first choice is simple and its numerical evaluation very
fast. However, for high polynomial degree $k$ the polynomial extrapolation,
possibly, becomes increasingly numerically unstable and the second choice
might be preferable for such cases.
\end{example}

Next, we prove optimal convergence rates with respect to the mesh width of the
postprocessed pressure $p_{M}^{\ast}$ from \eqref{defpostprocessedpressure}
for the functional $J_{\mathbf{z}}$ in Example \ref{ExJ}(1)
so that $f_{\mathbf{z}}$ is given by \eqref{Jzprimi}.

\begin{lemma}
Let Assumption \ref{AExt} hold and let the functionals $J_{\mathbf{z}}$ and
$f_{\mathbf{z}}$ be defined by \eqref{Jz0primi} and \eqref{Jzprimi}. Then,
Assumption \ref{Astruct} is satisfied.
\end{lemma}

%

\begin{proof}
Proposition \ref{Propcritfct}, \eqref{Jz0primi}, and an inverse inequality
(see, e.g., \cite{SchwabhpBook}) verify for any $q\in\mathbb{P}_{k-1}\left(
\mathbb{R}^{2}\right)  $ and $\mathbf{z}\in\mathcal{SC}_{\mathcal{T}}$ that
\[
\left\vert J_{\mathbf{z}}\left(  q\right)  \right\vert =\frac{\left\vert
q\left(  \mathbf{z}\right)  \right\vert }{\left\vert \left.  b_{k-1,\mathbf{z}%
}\right\vert _{K_{\mathbf{z}}}\left(  \mathbf{z}\right)  \right\vert }%
=\frac{\left\vert q\left(  \mathbf{z}\right)  \right\vert }{\binom{k+1}%
{2}\left\vert K_{\mathbf{z}}\right\vert ^{-1}}\leq2\frac{h_{\mathbf{z}}^{2}%
}{k^{2}}\left\Vert q\right\Vert _{L^{\infty}\left(  K_{\mathbf{z}}\right)
}\leq C_{\operatorname{inv}} h_{\mathbf{z}}\left\Vert q\right\Vert
_{L^{2}\left(  K_{\mathbf{z}}\right)  },
\]
where $C_{\operatorname{inv}} >0$ is independent of $h_{\mathbf{z}}$ and $k$.
Next, we use (\ref{bkL2}) to obtain
\[
\left\Vert b_{k-1,\mathbf{z}}\right\Vert _{L^{2}\left(  \omega_{\mathbf{z}%
}\right)  }^{2}=\sum_{K\in\mathcal{T}_{\mathbf{z}}}\left\vert K\right\vert
^{-1}\leq C_{\operatorname{sr}}^{2}h_{\mathbf{z}}^{-2}%
\]
for a constant $C_{\operatorname{sr}}$ depending only on the shape regularity
of the mesh. This implies the continuity,%
\[
\left\vert J_{\mathbf{z}}\left(  q\right)  \right\vert \leq\frac
{C_{\operatorname{inv}} C_{\operatorname{sr}}}{\left\Vert b_{k-1,\mathbf{z}%
}\right\Vert _{L^{2}\left(  \omega_{\mathbf{z}}\right)  }}\left\Vert
q\right\Vert _{L^{2}\left(  K_{\mathbf{z}}\right)  }.
\]
By setting $\mathcal{U}_{\mathbf{z}}:=K_{\mathbf{z}}$ we see that
(\ref{Eq:Jzcontestimate}) is satisfied for $C_{\mathbf{z}}:=C_{0}C_{1}$. We
choose $K_{\mathbf{z}}^{\operatorname*{ext}}:=K_{\mathbf{z}}$ and
$K_{\mathbf{z}}^{\prime,\operatorname*{ext}}$ as a minimal triangle which
contains $K_{\mathbf{z}}\cup K_{\mathbf{z}}^{\prime}$ so that (\ref{Kzext})
holds for a constant $\delta_{\mathbf{z}}=O\left(  1\right)  $ which only
depends on the shape regularity of the mesh. It remains to notice that the
definition (\ref{Jzprimi}) is of the form (\ref{formfz}).%
\end{proof}

A combination of Theorems \ref{Thm:ApproxPropertie pp-pressure} and
\ref{Thm:Recoverd bestapprox} implies optimal convergence rates in terms of the mesh width $h_{\mathcal{T}}$ of our simple postprocessing for the Scott-Vogelius element.

\begin{theorem}
Let Assumption \ref{AExt} hold and let the functionals $J_{\mathbf{z}}$ and
$f_{\mathbf{z}}$ be defined by (\ref{Jz0primi}) and (\ref{Jzprimi}). Let the
assumptions in Theorem \ref{Thm:ApproxPropertie pp-pressure} be satisfied.
Then there exists a constant $C>0$ depending on the shape-regularity of
$\mathcal{T}$, the constants $C_{\operatorname{Stein}}$,
$C_{\operatorname*{ov}}$, the domain $\Omega$, and the polynomial degree
$k\geq4$ such that the postprocessed pressure $p_{M}^{\ast}$ from
(\ref{defpostprocessedpressure}) satisfies%
\[
\left\Vert p-p_{M}^{\ast}\right\Vert _{L^{2}\left(  \Omega\right)  }\leq
Ch_{\mathcal{T}}^{\min\left\{  k,s-1\right\}  }\left\Vert p\right\Vert
_{H^{s-1}\left(  \Omega\right)  }.
\]

\end{theorem}

\section{Pressure-improvement for the pressure-wired Stokes element}
\label{SecPWStokes}
The pressure-wired Stokes element introduced in~\cite{Sauter_eta_wired} generalises
the standard Scott-Vogelius element by restricting the discrete
pressure space additionally at \textit{nearly singular vertices}
$\mathbf{z}\in\mathcal{V}\left(  \mathcal{T}\right)$ to guarantee a mesh-robust inf-sup stability.
In general, the divergence of the discrete
velocity field no longer vanishes pointwise while its marginality has
been analyzed in detail in \cite[Sec.~5]{Sauter_eta_wired}.
This section discusses a modified pressure space with a parameter $\eta\geq0$ introduced
below in full analogy to Section~\ref{Sec:Recover approx SV} for optimal approximation properties.
%

\begin{definition}
For $\eta\geq0$, the set of $\eta$\emph{-critical vertices} is given by%
\begin{equation}
\mathcal{C}_{\mathcal{T}}\left(  \eta\right)  :=\left\{  \left.  \mathbf{z}%
\in\mathcal{V}\left(  \mathcal{T}\right)  \;\right\vert \;\Theta\left(
\mathbf{z}\right)  \leq\eta\right\}  \label{Eq:Def eta critical vertices}%
\end{equation}
and the subset of $\eta$\emph{-super critical vertices} by%
\[
\mathcal{SC}_{\mathcal{T}}\left(  \eta\right)  :=\left\{  \left.
\mathbf{z}\in\mathcal{C}_{\mathcal{T}}\left(  \eta\right)  \;\right\vert
\;N_{\mathbf{z}}=1,3\right\}  .
\]

\end{definition}

Definition \ref{Def:Robinson vertex} generalises to $\eta$-super critical
vertices: a vertex $\mathbf{z\in}\mathcal{SC}_{\mathcal{T}}\left(
\eta\right)  $ is called a \textit{Robinson vertex} if (i): $\mathcal{T}%
_{\mathbf{z}}^{\ast}\cap\mathcal{T}_{\mathbf{y}}^{\ast}=\emptyset$ for all
$\mathbf{y}\in\mathcal{SC}_{\mathcal{T}}\left(  \eta\right)  \setminus\left\{
\mathbf{z}\right\}  $ and (ii): $\omega_{\mathbf{z}}^{\ast}\cap\mathcal{C}%
_{\mathcal{T}}\left(  \eta\right)  =\left\{  \mathbf{z}\right\}  $ hold.

\begin{remark}\label{rem:eta0}
For sufficiently small $\eta_{0}$ (depending only on the shape regularity of
the mesh and on $\Omega$) and $\eta\in\left[  0,\eta_{0}\right[  $ only the
four types of $\eta$-critical vertex configurations depicted in Figure
\ref{fig:Critical vertex_patch} are possible; see \cite[Lem. 2.13]%
{SAUTER202349} for details. Notice that $\mathcal{C}_{\mathcal{T}}\left(  0\right)  =\mathcal{C}_{\mathcal{T}}$ and
$\mathcal{SC}_{\mathcal{T}}\left(  0\right)  =\mathcal{SC}%
_{\mathcal{T}}$ hold.
Since $\mathcal{C}_{\mathcal{T}}(\eta) = \mathcal{C}_{\mathcal{T}}(1)$ for all $\eta\geq 1$ from $\Theta(\mathbf{z})\leq
1$, $\eta\leq 1$ in the following is not a restriction.
\end{remark}

\begin{definition}
	[pressure-wired Stokes element, {\cite[Lem.~1]{Sauter_eta_wired}}]
Given $k\in\mathbb N$ and $0\leq \eta\leq 1$, 
the \emph{pressure-wired Stokes element} $\left(  \mathbf{S}_{k,0}\left(  \mathcal{T}\right)  ,M_{\eta,k-1}\left(
\mathcal{T}\right)  \right)$  is defined with the \emph{reduced pressure space} 
\begin{equation}
M_{\eta,k-1}\left(  \mathcal{T}\right)  :=\left\{  \left.  q\in\mathbb{P}%
_{k-1,0}\left(  \mathcal{T}\right)  \;\right\vert \;\forall\mathbf{z}%
\in\mathcal{C}_{\mathcal{T}}\left(  \eta\right)  :\ A_{\mathcal{T},\mathbf{z}%
}\left(  q\right)  =0\right\}  . \label{Eq:Def Meta}%
\end{equation}

\end{definition}

As already announced in the introduction there is a significant loss in accuracy of the pressure
approximation for the pressure-wired Stokes element if the mesh contains
$\eta$-super critical vertices. As a remedy we introduce a modification in
complete analogy to Section~\ref{Sec:Recover approx SV} for the Scott-Vogelius element.

\begin{definition}
[modified pressure-wired Stokes element]
Given $k\in\mathbb N, 0\leq \eta\leq 1$, and
functionals 
$f_{\mathbf{z}}:\mathbb{P}_{k-1}\left(  \mathcal{T}\right)  \rightarrow\mathbb{R}$ with~\eqref{estfz} for all 
$\mathbf{z}\in\mathcal{SC}%
_{\mathcal{T}}\left(  \eta\right)  $, the \emph{modified
pressure-wired Stokes element} $\left(  \mathbf{S}_{k,0}\left(  \mathcal{T}\right)  ,M_{\eta,k-1}%
^{\operatorname{mod}}\left(  \mathcal{T}\right)  \right)$ is given by the modified reduced pressure space 
\begin{equation}
M_{\eta,k-1}^{\operatorname{mod}}\left(  \mathcal{T}\right)  :=\left\{
q+\sum_{\mathbf{z}\in\mathcal{SC}_{\mathcal{T}}\left(  \eta\right)
}f_{\mathbf{z}}\left(  q\right)  \left(  b_{k-1,\mathbf{z}}\right)
_{\operatorname*{mvz}}\;:\;q\in M_{\eta,k-1}\left(  \mathcal{T}\right)
\right\}  . \label{DefModPressSp}%
\end{equation}

\end{definition}

\subsection{Stability of the modified pressure-wired Stokes
element\label{SecStabmodPWStokes}}

The modified pressure-wired Stokes element inherits the mesh-robust discrete inf-sup stability.

\begin{figure}[ptb]
\centering
\resizebox{\textwidth}{8em}{
	\begin{tikzpicture}[scale=0.8]
				\draw (0,0) node[] (z) {$\mathbf{z}$};
		\coordinate (A) at (4,0);
		\coordinate (B) at ($(z)!0.5!80:(A)$);
		\coordinate (C) at ($(z)!0.9!200:(A)$);
		\coordinate (D) at ($(z)!0.53!300:(A)$);
		\draw (A) -- (B) -- (C) -- (D) -- cycle;
		\draw (z) -- (A) node[midway,above] {$E_1$};
		\draw (z) -- (B) node[midway,right] {$E_2$};
		\draw (z) -- (C) node[midway,above] {$E_3$};
		\draw (z) -- (D) node[midway,left] {$E_4$};
		\node (T1) at (barycentric cs:z=1,A=1,B=1) {$K_1$};
		\node (T1) at (barycentric cs:z=1,C=1,B=1) {$K_2$};
		\node (T1) at (barycentric cs:z=1,A=1,D=1) {$K_4$};
		\node (T1) at (barycentric cs:z=1,C=1,D=1) {$K_3$};
	\end{tikzpicture}
	\begin{tikzpicture}[scale=0.8]
				\draw (0,0) node[] (z) {$\mathbf{z}$};
		\coordinate (A) at (4,0);
		\coordinate (B) at ($(z)!0.5!80:(A)$);
		\coordinate (C) at ($(z)!0.9!200:(A)$);
		\coordinate (D) at ($(z)!0.53!300:(A)$);
		\draw (A) -- (B) -- (C) -- (D);
		\draw (z) -- (A) node[midway,above] {$E_1$};
		\draw (z) -- (B) node[midway,right] {$E_2$};
		\draw (z) -- (C) node[midway,above] {$E_3$};
		\draw (z) -- (D) node[midway,left] {$E_4$};
		\node (T1) at (barycentric cs:z=1,A=1,B=1) {$K_1$};
		\node (T1) at (barycentric cs:z=1,C=1,B=1) {$K_2$};
				\node (T1) at (barycentric cs:z=1,C=1,D=1) {$K_3$};
		\fill[pattern=north east lines, ] ([xshift=0.35em]z.south) -- (D) -- ([xshift=.3em]D) -- ([xshift=.65em]z.south);
		\fill[pattern=north east lines, ] (z.east) rectangle ([yshift=-.3em]A);
	\end{tikzpicture}
	\begin{tikzpicture}[scale=0.8]
				\draw (0,0) node[] (z) {$\mathbf{z}$};
		\coordinate (A) at (4,0);
		\coordinate (B) at ($(z)!0.5!80:(A)$);
		\coordinate (C) at ($(z)!0.9!200:(A)$);
		\coordinate (D) at ($(z)!0.53!300:(A)$);
		\draw (A) -- (B) -- (C);		\draw (z) -- (A) node[midway,above] {$E_1$};
		\draw (z) -- (B) node[midway,right] {$E_2$};
		\draw (z) -- (C) node[midway,above] {$E_3$};
				\node (T1) at (barycentric cs:z=1,A=1,B=1) {$K_1$};
		\node (T1) at (barycentric cs:z=1,C=1,B=1) {$K_2$};
		\phantom{		\node (T1) at (barycentric cs:z=1,A=1,D=1) {$K_4$};
		\node (T1) at (barycentric cs:z=1,C=1,D=1) {$K_3$};}
		\fill[pattern=north west lines, ] ([yshift=-0.25em]z.west) -- (C) -- ([yshift=-.3em]C) -- ([yshift=-.55em]z.west);
		\fill[pattern=north east lines, ] (z.east) rectangle ([yshift=-.3em]A);
	\end{tikzpicture}
	\begin{tikzpicture}[scale=0.8]
				\draw (0,0) node[] (z) {$\mathbf{z}$};
		\coordinate (A) at (4,0);
		\coordinate (B) at ($(z)!0.5!80:(A)$);
		\coordinate (C) at ($(z)!0.9!200:(A)$);
		\coordinate (D) at ($(z)!0.53!300:(A)$);
		\draw (A) -- (B);		\draw (z) -- (A) node[midway,above] {$E_1$};
		\draw (z) -- (B) node[midway,right] {$E_2$};
						\phantom{		\node (T1) at (barycentric cs:z=1,A=1,D=1) {$K_4$};
		\node (T1) at (barycentric cs:z=1,C=1,D=1) {$K_3$};}
		\node (T1) at (barycentric cs:z=1,A=1,B=1) {$K_1$};
		\fill[pattern=north east lines, ] (z.east) rectangle ([yshift=-.3em]A);
		\fill[pattern=north west lines, ] ([xshift=0.08em]z.north) -- (B) -- ([xshift=-.3em]B) -- ([xshift=-.22em]z.north);
						\end{tikzpicture}}
						\caption{Vertex patch for an interior $\eta$-critical
vertex $\mathbf{z}\in\mathcal{V}_{\Omega}(\mathcal{T})$ with $N_{\mathbf{z}%
}=4$ (resp.~boundery $\eta$-critica vertex $\mathbf{z}\in\mathcal{V}%
_{\partial\Omega}(\mathcal{T})$ with $N_{\mathbf{z}}=1,2,3$) triangles}%
\label{fig:Critical vertex_patch}%
\end{figure}

\begin{lemma}
\label{Prop:Stability pwS mod} 
The modified pressure-wired Stokes element 
is inf-sup stable for all $k\geq 4$ and $0\leq \eta\le 1$ with $\beta\left(  \mathbf{S} _{k,0}\left(  \mathcal{T}\right)  ,M_{\eta,k-1}^{\operatorname{mod}}\left(
\mathcal{T}\right)  \right)  \geq c(\Theta_{\min}+\eta)/C_{f}$ for
$C_f=1+\sum^{}_{\mathbf{z}\in\mathcal{S}\mathcal{C}_{\mathcal{T}}(\eta)} C_{f_{\mathbf{z}}}$. The constant $c>0$
exclusively depends on the shape regularity of the mesh.
\end{lemma}

\begin{proof}
	The proof of Lemma~\ref{Prop:Stability pwS mod} is a modification of that of Lemma~\ref{Propinfsupmod} and given for
	completeness.
	The surjective extrapolation operator $\mathcal{E}_{\eta,k}:M_{\eta,k-1}\left(  \mathcal{T}\right)\to
	M_{\eta,k-1}^{\operatorname*{mod}}\left(  \mathcal{T}\right)  $ given by
\[
\mathcal{E}_{\eta,k}q=q+\sum_{\mathbf{z}\in\mathcal{SC}_{\mathcal{T}}\left(
\mathcal{\eta}\right)  }f_{\mathbf{z}}\left(  q\right)  \left(
b_{k-1,\mathbf{z}}\right)  _{\operatorname*{mvz}}\qquad\forall q\in M_{\eta,k-1}\left(
	\mathcal{T}\right)
\]
generalises
$\mathcal{E}_{k}$ from \eqref{defextraop}.
Lemma~2 in~\cite{Sauter_eta_wired} provides, for any $q\in M_{\eta,k-1}\left(\mathcal{T}\right)$,
some $\mathbf{v} _{q}:=\Pi_{k}q\in \mathbf{S}_{k,0}(\mathcal{T})$ 
with $\operatorname{div}(\mathbf{v}_q)=q$ and
$c(\Theta_{\min}+\eta)\|\mathbf{v}_q\|_{\mathbf{H}^1(\Omega)}\leq \|q\|_{L^2(\Omega)}$.
As in the proof of Lemma \ref{Propinfsupmod},
\[
\left(  \mathcal{E}_{\eta,k}q,q\right)
_{L^{2}\left(  \Omega\right)  }=\|q\|_{L^2(\Omega)}^2\quad\text{and\quad}\left\Vert \mathcal{E}_{\eta
,k}q\right\Vert _{L^{2}\left(  \Omega\right)  }\leq C_f  \left\Vert q\right\Vert _{L^{2}\left(
\Omega\right)  } \qquad\forall q\in M_{\eta,k-1}\left(\mathcal{T}\right),
\]
the right-inverse of the divergence, and the surjectivity of $\mathcal{E}_{\eta,k}$ bound the inf-sup constant by
\begin{align*}
\inf_{q^{\ast}\in M_{\eta,k-1}^{\operatorname{mod}}\left(  \mathcal{T}\right)
\setminus\left\{  0\right\}  }\sup_{\mathbf{v}\in\mathbf{S}_{k,0}\left(
\mathcal{T}\right)  \setminus\left\{  \mathbf{0}\right\}  }\frac{\left(
q^{\ast},\operatorname{div}\mathbf{v}\right)  _{L^{2}\left(  \Omega\right)  }%
}{\left\Vert q^{\ast}\right\Vert _{L^{2}\left(  \Omega\right)  }\left\Vert
\mathbf{v}\right\Vert _{H^{1}\left(  \Omega\right)  }}  &  \geq\inf_{q\in
M_{\eta,k-1}\left(  \mathcal{T}\right)  \setminus\left\{  0\right\}  }%
\frac{\left(  \mathcal{E}_{\eta,k}q,q\right)
_{L^{2}\left(  \Omega\right)  }}{\left\Vert \mathcal{E}_{\eta,k}q\right\Vert
_{L^{2}\left(  \Omega\right)  }\left\Vert \mathbf{v}_{q}\right\Vert
_{H^{1}\left(  \Omega\right)  }}\\
&  \geq c(\Theta_{\min}+\eta)/C_f .\qedhere
\end{align*}
\end{proof}
A consequence of Lemma~\ref{Prop:Stability pwS mod} and the theory of mixed methods~\cite[Chap.~5, Thm.~5.2.3]{BoffiBrezziFortin} is the
quasi-optimality of the discrete solution. Recall $C_f$ from Lemma~\ref{Prop:Stability pwS mod}.

\begin{corollary}
\label{Cor:Best approx pwS mod} Given $k\geq4$ and $0\leq \eta\leq 1$, let $\left(
\mathbf{u},p\right)  \in\mathbf{H}_{0}^{1}\left(  \Omega\right)  \times
L_{0}^{2}\left(  \Omega\right)  $ solve
\eqref{varproblemstokes}. The discrete solution $\left(  \mathbf{u}_{S},p_{M}^{\ast}\right)$ of \eqref{discrStokes} for the
choice $\left(  \mathbf{S},M\right)  =\left(  \mathbf{S}_{k,0}\left(
\mathcal{T}\right)  ,M_{\eta,k-1}^{\operatorname{mod}}\left(  \mathcal{T}%
\right)  \right)  $
satisfies the quasi-optimal error estimate
\begin{subequations}
\label{Corests}
\end{subequations}%
\begin{align}
\left\Vert \nabla \left(\mathbf{u}-\mathbf{u}_{\mathbf{S}} \right)\right\Vert _{\mathbb{L}^{2}\left(  \Omega\right)  }&\leq   C\left(  \frac{C_f
}{\Theta_{\min}+\eta}\inf_{\mathbf{v}\in\mathbf{S}%
_{k,0}\left(  \mathcal{T}\right)  }\left\Vert \mathbf{u}-\mathbf{v}\right\Vert
_{\mathbf{H}^{1}\left(  \Omega\right)  }+\inf_{q\in M_{\eta,k-1}^{\operatorname{mod}%
}\left(  \mathcal{T}\right)  }\left\Vert p-q\right\Vert _{L^{2}\left(
\Omega\right)  }\right)  ,\tag{%
\ref{Corests}%
a}\label{Corestsa}\\
\left\Vert p-p_{M}^{\ast}\right\Vert _{L^{2}\left(  \Omega\right)  }&\leq 
\frac{C_{\operatorname{vel}}C_f^2}%
{\left(  \Theta_{\min}+\eta\right)  ^{2}}\inf_{\mathbf{v}\in
\mathbf{S}_{k,0}\left(  \mathcal{T}\right)  }\left\Vert \mathbf{u}%
-\mathbf{v}\right\Vert _{\mathbf{H}^{1}\left(  \Omega\right)  }+\tag{%
\ref{Corests}%
b}\label{Corestsb}
\frac{C_{\operatorname{pres}}C_f}{\Theta_{\min
}+\eta}\inf_{q\in M_{\eta,k-1}^{\operatorname{mod}}\left(
\mathcal{T}\right)  }\left\Vert p-q\right\Vert _{L^{2}\left(  \Omega\right)
}.\nonumber
\end{align}
The positive
constants $C,C_{\operatorname{vel}},C_{\operatorname{pres}}$ only depend on
the shape-regularity of the mesh and the domain $\Omega$.\qed
\end{corollary}

\subsection{Optimal convergence rates\label{SecPWStokesConv}}

Convergence rates for the discrete solution $\left(  \mathbf{u}_{\mathbf{S}}%
,p_{M}^{\ast}\right)  $ follow from Corollary \ref{Cor:Best approx pwS mod}
once the approximation property of $M_{\eta,k-1}^{\operatorname{mod}}\left(
\mathcal{T}\right)$ is clarified.
One key argument is an analog characterisation of the orthogonal complement of $M_{\eta,k-1}(\mathcal{T})$ in $\mathbb
P_{k-1,0}(\mathcal{T})$ as in Proposition~\ref{PropOrthCompl}.

\begin{lemma}\label{PropOrthComplEta}
Let all $\eta$-super critical
vertices be of Robinson type. Then the decomposition
\begin{align*}
	\mathbb P_{k-1,0}(\mathcal{T}) = M_{\eta,k-1}(\mathcal{T})\oplus \operatorname{span}\{b_{k-1,\mathbf{z}}\ |\ \mathbf{z}\in
	\mathcal{C}_{\mathcal{T}}(\eta)\setminus\mathcal{S}\mathcal{C}_{\mathcal{T}}(\eta)\}\oplus
		\operatorname{span}\{(b_{k-1,\mathbf{y}})_{\operatorname{mvz}}\ |\
		\mathbf{y}\in\mathcal{S}\mathcal{C}_{\mathcal{T}}(\eta)\}
\end{align*}
is $L^2$ orthogonal, i.e., any $\mathbf{z}\in
\mathcal{C}_{\mathcal{T}}(\eta)\setminus\mathcal{S}\mathcal{C}_{\mathcal{T}}(\eta)$, $\mathbf{y}\in
\mathcal{S}\mathcal{C}_{\mathcal{T}}(\eta)$, and $q_M\in M_{\eta,k-1}(\mathcal{T})$ satisfy
\begin{align*}
	(q_M, b_{k-1,\mathbf{z}})_{L^2(\Omega)}=(q_M, \left(  b_{k-1,\mathbf{y}}\right)
_{\operatorname*{mvz}})_{L^2(\Omega)}=(b_{k-1,\mathbf{z}}, \left(  b_{k-1,\mathbf{y}}\right)
_{\operatorname*{mvz}})_{L^2(\Omega)}=0.
\end{align*}
\end{lemma}
\begin{proof}
	Since $N_{\mathbf{z}}$ is even for all
	$\mathbf{z}\in\mathcal{C}_{\mathcal{T}}(\eta)\setminus\mathcal{S}\mathcal{C}_{\mathcal{T}}(\eta)$ by definition, 
	\eqref{bkintegral} in Proposition~\ref{Propcritfct} verifies $\overline{b_{k-1,\mathbf{z}}}=0$ so that $b_{k-1,\mathbf{z}}=\left(  b_{k-1,\mathbf{z}}\right)
_{\operatorname*{mvz}}$.
This, the orthogonality $(q_M, \left(  b_{k-1,\mathbf{y}}\right)
_{\operatorname*{mvz}})_{L^2(\Omega)}$ for all $q_M\in M_{\eta,k-1}(\mathcal{T})$ and
$\mathbf{z}\in\mathcal{C}_{\mathcal{T}}(\eta)$ as in Proposition~\ref{Propcritfct}, and Remark~\ref{Rem:Avoiding
techinacl difficulties}(i) adapted to $\eta$-critical vertices conclude the proof; further details are omitted.
\end{proof}
As in Section \ref{Sec:Recover approx SV}, we employ a
representation of the space $\widetilde{Q}_{0}^{h,k}$ with optimal approximation properties from~\eqref{defQtilde} in terms of functions
in $M_{\eta,k-1}\left(  \mathcal{T}\right)  $ and $\left\{  \left.
b_{k-1,\mathbf{z}}\;\right\vert \;\mathbf{z}\in\mathcal{C}_{\mathcal{T}%
}\left(  \eta\right)  \right\}  $.
\begin{lemma}
\label{Lem:Decomposition of Q0hk wrt Meta}Let all $\eta$-super critical
vertices be of Robinson type and let $k\geq4$ be given. Then for all
$\widetilde{q} \in\widetilde{Q}_{0}^{h,k}$ there exist $q_{0}\in M_{\eta,k-1}\left(
\mathcal{T}\right)  $ and $\theta_{\mathbf{z}}\in\mathbb{R}$ for all
$\mathbf{z}\in\mathcal{SC}_{\mathcal{T}}\left(  \eta\right)  $ such that
\begin{equation} \label{Eq:Decomposition of Q0hk wrt Meta}
\widetilde{q} =q_{0}+\sum_{\mathbf{z}\in\mathcal{SC}_{\mathcal{T}}\left(  \eta\right)
}\theta_{\mathbf{z}}\left(  b_{k-1,\mathbf{z}}\right)  _{\operatorname*{mvz}%
}.
\end{equation}

\end{lemma}%

\begin{proof}
	The proof is a straightforward modification of that of Proposition~\ref{Lem:Decomposition of Q0hk} with the orthogonal
	decomposition of $\mathbb P_{k-1,0}(\mathcal{T})$ from Lemma~\ref{PropOrthComplEta}; further details are omitted.
\end{proof}

\begin{theorem}
\label{Thm:Recoverd approx Meta mod}Let all $\eta$-super critical vertices be
of Robinson type and let Assumption \ref{Astruct} be satisfied. For any $p\in
H^{s-1}\left(  \Omega\right)  \cap L_{0}^{2}\left(  \Omega\right)  $ with
$s>1$, there exists $p_{M}^{\ast}\in M_{\eta,k-1}^{\operatorname*{mod}%
}\left(  \mathcal{T}\right)  $ such that%
\begin{equation}
\left\Vert p-p_{M}^{\ast}\right\Vert _{L^{2}\left(  \Omega\right)  }\leq
C_{\operatorname*{apx}}\frac{\left(  \left(  1+2\delta_{\max}\right)
h_{\mathcal{T}}\right)  ^{\min\left\{  k,s-1\right\}  }}{k^{s-1}}\left\Vert
p\right\Vert _{H^{s-1}\left(  \Omega\right)  } \label{errorestmod}%
\end{equation}
for $C_{\operatorname*{apx}}:=C\sqrt{\operatorname{card}\mathcal{SC}_{\mathcal{T}%
}\left(  \eta\right)  }C_{\operatorname{Stein}}C_{J}T_{k-1}\left(
1+c\delta_{\max}\right)$
with $C_{J}$ from
\eqref{CJconstant} and $\delta_{\max}$ from \eqref{CJconstant}.
The constant $C$ depends only on $C_{\operatorname*{ov}}$ from~\eqref{defCo}
and on the shape regularity of the mesh.
\end{theorem}

The proof of this theorem follows by applying the arguments in
the proof of Theorem \ref{Thm:Recoverd bestapprox} verbatim to $\mathcal{SC}%
_{\mathcal{T}}\left(  \eta\right)  $ and Lemma
\ref{Lem:Decomposition of Q0hk wrt Meta} instead of $\mathcal{SC}%
_{\mathcal{T}}$ and Lemma \ref{Lem:Decomposition of Q0hk}.

\subsection{Control of $\operatorname*{div}\mathbf{u}_{\mathbf{S}}%
$\label{DivControl}}

The pressure-wired Stokes element reduces the pressure space not only at exactly
singular vertices $\mathcal{C}_{\mathcal{T}}$ but also at $\eta
$-critical vertices $\mathbf{z}\in\mathcal{C}_{\mathcal{T}}\left(
\eta\right)  \backslash\mathcal{C}_{\mathcal{T}}$. Therefore, the discrete
velocity field $\mathbf{u}_{\mathbf{S}}\in\mathbf{S}_{k,0}\left(
\mathcal{T}\right)  $ of the pressure-wired Stokes element is in general not
pointwise divergence free. 
However, Theorem 3 in \cite{Sauter_eta_wired}
guarantees that $\left\Vert \operatorname{div}\mathbf{u}_{\mathbf{S}%
}\right\Vert _{L^{2}\left(  \Omega\right)  }$ tends to zero at least linearly
in $\eta$. An analogous statement holds for the modified
pressure-wired Stokes element. Consider the open subset 
\begin{align}\label{eqn:Omega_star_def}
	\Omega(\eta):= \bigcup_{\mathbf{z}\in\mathcal{C}_{\mathcal{T}}(\eta)}\operatorname{int}(\omega_\mathbf{z})\subset\Omega
\end{align}
for $0\leq \eta\leq 1$, where $\operatorname{int}(\omega_{\mathbf{z}})$ denotes the interior of the vertex patch $\omega_{\mathbf{z}}$
from~\eqref{nodalpatch}, and define
\begin{equation}
\mathbf{S}_{\eta,k,0}\left(  \mathcal{T}\right)  :=\left\{  \mathbf{v}%
\in\mathbf{S}_{k,0}\left(  \mathcal{T}\right)  \mid A_{\mathcal{T},\mathbf{z}%
}\left(  \operatorname{div}\mathbf{v}\right)  =0\quad\forall\mathbf{z}%
\in\mathcal{C}_{\mathcal{T}}\left(  \eta\right)  \right\}  .
\label{DefSetak0T}%
\end{equation}

\begin{theorem}
	\label{Thm:Divergence estimate}
	There exists $\eta_0>0$ such that the following holds: Given $0\leq\eta < \eta_0$, let Assumption \ref{Astruct} be satisfied for all
	$\eta$-super critical vertices $\mathbf{z}\in\mathcal{S}\mathcal{C}_{\mathcal{T}}(\eta)$
	that are additionally assumed to be of Robinson type.
For $k\geq4 $, the discrete solution $\left(  \mathbf{u}_{\mathbf{S}},p_{M}^{\ast}\right)
\in\left(  \mathbf{S}_{k,0}\left(  \mathcal{T}\right)  ,M_{\eta,k-1}%
^{\operatorname{mod}}\left(  \mathcal{T}\right)  \right)  $ to
(\ref{discrStokes}) satisfies
\begin{equation}
\left\Vert \operatorname{div}\mathbf{u}_{\mathbf{S}}\right\Vert _{L^{2}\left(
\Omega\right)  }\leq C_{\operatorname*{div}}\eta\inf_{\mathbf{w}_{\mathbf{S}%
}\in\mathbf{S}_{\eta,k,0}\left(  \mathcal{T}\right)  }\left\Vert \nabla\left(
\mathbf{u}_{\mathbf{S}}-\mathbf{w}_{\mathbf{S}}\right)  \right\Vert
_{\mathbb{L}^{2}\left(  \Omega(\eta)\right)  }. \label{divest}%
\end{equation}
The constant $C_{\mathrm{div}}>0$ in (\ref{divest}) only depends on the shape-regularity
constant and the domain $\Omega$.
\end{theorem}
The remaining parts of this section are devoted to the proof of Theorem~\ref{Thm:Divergence estimate}.
Since $\mathbf{u}_{\mathbf{S}}$ solves \eqref{discrStokes} and has integral mean zero, its divergence is orthogonal to ${M_{\eta,k-1}^{\operatorname{mod}%
}\left(  \mathcal{T}\right)  }$ in $P_{k-1,0}(\mathcal{T})$, namely
\begin{align}\label{eqn:div_pw_orthogonality}
	\operatorname{div}\mathbf{u}_{\mathbf{S}}\in{M_{\eta,k-1}^{\operatorname{mod}%
			}\left(  \mathcal{T}\right)  }^{\perp}:=\left\{  \left.  q\in\mathbb{P}%
				_{k-1,0}\left(  \mathcal{T}\right)  \;\right\vert \;\forall p\in M_{\eta
				,k-1}^{\operatorname{mod}}\left(  \mathcal{T}\right)  :\ \left(  q,p\right)
			_{L^{2}\left(  \Omega\right)  }=0\right\}.
\end{align}
We characterise the
$L^{2}\left(  \Omega\right)  $-orthogonal
complement $M_{\eta,k-1}^{\operatorname{mod}}\left(  \mathcal{T}\right)
^{\perp}$ to prove
Theorem~\ref{Thm:Divergence estimate}. Let $\phi_{\mathbf{z}}\in M_{\eta
,k-1}\left(  \mathcal{T}\right)  $ be the Riesz representative of 
$f_{\mathbf{z}}:M_{\eta,k-1}\left(  \mathcal{T}\right)
\rightarrow\mathbb{R}$, i.e., $\phi_{\mathbf{z}}$ satisfies
\begin{equation}
\left(  \phi_{\mathbf{z}},q_{0}\right)  _{L^{2}\left(  \Omega\right)
}=f_{\mathbf{z}}\left(  q_{0}\right)  \qquad\forall q_{0}\in M_{\eta
,k-1}\left(  \mathcal{T}\right)  . \label{Eq:Riesz rep}%
\end{equation}

\begin{lemma}
\label{Prop:Orhtogonal complement Mmod}
Under the assumptions of Theorem
\ref{Thm:Divergence estimate}, a basis for $M_{\eta,k-1}^{\operatorname{mod}}\left(  \mathcal{T}\right)  ^{\perp
}=\operatorname{span}\left(  \mathcal{B}\right)$ is given by $\mathcal{B}:=\left\{  \left.
b_{k-1,\mathbf{z}}\;\right\vert \;\forall\mathbf{z}\in\mathcal{C}%
_{\mathcal{T}}\left(  \eta\right)  \setminus\mathcal{SC}_{\mathcal{T}}\left(
\eta\right)  \right\}  \cup\left\{  \left.  \varphi_{\mathbf{z}}\;\right\vert
\;\forall\mathbf{z}\in\mathcal{SC}_{\mathcal{T}}\left(  \eta\right)  \right\}
$ with
\begin{equation}
\label{Eq:Def varphi}\varphi_{\mathbf{z}}:=\phi_{\mathbf{z}}-\frac
{1}{\left\Vert b_{k-1,\mathbf{z}}\right\Vert _{L^{2}\left(  \Omega\right)
}^{2}}\left(  \left(  b_{k-1,\mathbf{z}}\right)  _{\operatorname*{mvz}}%
+\sum_{\mathbf{y}\in\mathcal{SC}_{\mathcal{T}}\left(  \eta\right)  }\left(
\overline{b_{k-1,\mathbf{z}}},b_{k-1,\mathbf{y}}\right)  _{L^{2}\left(
\Omega\right)  }\phi_{\mathbf{y}}\right)  .
\end{equation}

\end{lemma}

%

\begin{proof}
A counting argument with Propostion \ref{Propcritfct} and Lemma \ref{PropOrthComplEta} show 
\begin{align*}
	\operatorname{card}
\mathcal{B} = \operatorname{card} \mathcal{C}_{\mathcal{T}} \left( \eta \right)=\operatorname{dim}
{M_{\eta,k-1}\left(  \mathcal{T}\right)  }^{\perp},
\end{align*}
where the orthogonal complement ${M_{\eta,k-1}\left(  \mathcal{T}\right)  }^{\perp}$ of ${M_{\eta,k-1}\left(
\mathcal{T}\right)  }$ in $\mathbb P_{k-1,0}(\mathcal{T})$ is defined analogously to~\eqref{eqn:div_pw_orthogonality}.
The definition of $M_{\eta
,k-1}^{\operatorname{mod}}\left(  \mathcal{T}\right)  $ in~\eqref{DefModPressSp} and
$\left(b_{k-1,\mathbf{z}}\right)_{\operatorname{mvz}} \in {M_{\eta,k-1} \left( \mathcal{T} \right)}^{\perp}$ for all
$\mathbf{z}\in\mathcal{SC}_{\mathcal{T}}(\eta)$
from Lemma~\ref{PropOrthComplEta} provide $\operatorname{dim} M_{\eta,k-1} \left( \mathcal{T} \right) = \operatorname{dim} M_{\eta,k-1}^{\operatorname{mod}} \left( \mathcal{T} \right)$.
Hence, $\operatorname{card} \mathcal{B} = \operatorname{dim} M_{\eta,k-1}^{\operatorname{mod}} \left( \mathcal{T} \right)^{\perp}$. 
The set $\left\{  \left.
b_{k-1,\mathbf{z}}\;\right\vert \;\forall\mathbf{z}\in\mathcal{C}%
_{\mathcal{T}}\left(  \eta\right)  \setminus\mathcal{SC}_{\mathcal{T}}\left(
\eta\right)  \right\}$ 
is linearly independent by Proposition~\ref{Propcritfct} and orthogonal to both, 
$M_{\eta
,k-1}^{\operatorname{mod}}\left(  \mathcal{T}\right)
\subset M_{\eta,k-1}(\mathcal{T})+\{\left(  b_{k-1,\mathbf{z}}\right)
_{\operatorname*{mvz}}\ |\ \mathbf{z}\in\mathcal{S}\mathcal{C}_{\mathcal{T}}(\eta)\}$
and $\mathcal{B}_0\coloneqq\{\varphi_{\mathbf{z}}\ |\
\forall\mathbf{z}\in\mathcal{S}\mathcal{C}_{\mathcal{T}}(\eta)\}\subset \mathcal{B}$,
by Lemma~\ref{PropOrthComplEta}.
To analyse the remaining subset
$\mathcal{B}_0$,
let
\[
q=q_{0}+\sum_{\mathbf{y}\in\mathcal{SC}_{\mathcal{T}}\left(  \eta\right)
}f_{\mathbf{y}}\left(  q_{0}\right)  \left(  b_{k-1,\mathbf{y}}\right)
_{\operatorname*{mvz}}\in M_{\eta,k-1}^{\operatorname{mod}} \left( \mathcal{T} \right)
\]
be arbitrary with
$q_{0}\in M_{\eta,k-1}\left(
\mathcal{T}\right)$.
Given any $\mathbf{z}\in\mathcal{SC}_{\mathcal{T}}\left(
\eta\right)  $, consider
\begin{align}\label{eqn:vfz_q_split}
\left(  \varphi_{\mathbf{z}},q\right)  _{L^{2}\left(  \Omega\right)  }=\left(
\varphi_{\mathbf{z}},q_{0}\right)  _{L^{2}\left(  \Omega\right)  }%
+\sum_{\mathbf{y}\in\mathcal{SC}_{\mathcal{T}}\left(  \eta\right)
}f_{\mathbf{y}}\left(  q_{0}\right)  \left(  \varphi_{\mathbf{z}},\left(
b_{k-1,\mathbf{y}}\right)  _{\operatorname*{mvz}}\right)  _{L^{2}\left(
\Omega\right)  }.
\end{align}
The definition of $\varphi_{\mathbf{z}}$ in \eqref{Eq:Def varphi}
with \eqref{Eq:Riesz rep} and $\left(  \left(
b_{\mathbf{z},k-1}\right)  _{\operatorname{mvz}},q_{0}\right)  _{L^{2}\left(
\Omega\right)  }=0$ by Lemma~\ref{PropOrthComplEta} reveal
\begin{align}
\left(
\varphi_{\mathbf{z}},q_{0}\right)  _{L^{2}\left(  \Omega\right)  }
&  =f_{\mathbf{z}}\left(  q_{0}\right)  -\frac{1}{\left\Vert b_{k-1,\mathbf{z}%
}\right\Vert _{L^{2}\left(  \Omega\right)  }^{2}}\sum_{\mathbf{y}%
\in\mathcal{SC}_{\mathcal{T}}\left(  \eta\right)  }\left(  \overline
{b_{k-1,\mathbf{z}}},b_{k-1,\mathbf{y}}\right)  _{L^{2}\left(  \Omega\right)
}f_{\mathbf{y}}\left(  q_{0}\right)  .\label{SIfinal}%
\end{align}
A similar computation with $\left(  \phi_{\mathbf{a}},\left(  b_{k-1,\mathbf{y}}\right)
_{\operatorname*{mvz}}\right)  _{L^{2}\left(  \Omega\right)
} = 0$ for all $\mathbf{a}\in\mathcal{S}\mathcal{C}_{\mathcal{T}}(\eta)$ from 
$\phi_{\mathbf{z}}\in M_{\eta,k-1}\left(  \mathcal{T}\right)$
and Lemma~\ref{PropOrthComplEta} implies for all $\mathbf{z}\in\mathcal{S}\mathcal{C}_{\mathcal{T}}(\eta)$ that
\begin{align*}
\left(  \varphi_{\mathbf{z}},\left(  b_{k-1,\mathbf{y}}\right)
_{\operatorname*{mvz}}\right)  _{L^{2}\left(  \Omega\right)  }
=-\frac{\left(  \left(
b_{k-1,\mathbf{z}}\right)  _{\operatorname*{mvz}},\left(  b_{k-1,\mathbf{y}%
}\right)  _{\operatorname*{mvz}}\right)  _{L^{2}\left(  \Omega\right)  }%
}{\left\Vert b_{k-1,\mathbf{z}}\right\Vert _{L^{2}\left(  \Omega\right)  }%
^{2}}
= &  -\delta_{\mathbf{z},\mathbf{y}}+\frac{\left(  \overline{b_{k-1,\mathbf{z}%
}},b_{k-1,\mathbf{y}}\right)  _{L^{2}\left(  \Omega\right)  }}{\left\Vert
b_{k-1,\mathbf{z}}\right\Vert _{L^{2}\left(  \Omega\right)  }^{2}}
\end{align*}
using the integral mean zero property of $\left(
b_{k-1,\mathbf{z}}\right)  _{\operatorname*{mvz}}$ and Remark~\ref{Rem:Avoiding techinacl difficulties}(i) as
all super-critical vertices are of Robinson type in the last step.
Here, $\delta_{\mathbf{z},\mathbf{y}}$ denotes the Kronecker delta.
The sum over all $\mathbf{y}\in\mathcal{S}\mathcal{C}_{\mathcal{T}}(\eta)$ results in
\begin{align}
\sum_{\mathbf{y}\in\mathcal{SC}_{\mathcal{T}}\left(  \eta\right)
}f_{\mathbf{y}}\left(  q_{0}\right)  \left(  \varphi_{\mathbf{z}},\left(
b_{k-1,\mathbf{y}}\right)  _{\operatorname*{mvz}}\right)  _{L^{2}\left(
\Omega\right)  }
  =-f_{\mathbf{z}}\left(  q_{0}\right)  +\sum_{\mathbf{y}\in\mathcal{SC}%
_{\mathcal{T}}\left(  \eta\right)  }f_{\mathbf{y}}\left(  q_{0}\right)
\frac{\left(  \overline{b_{k-1,\mathbf{z}}},b_{k-1,\mathbf{y}}\right)
_{L^{2}\left(  \Omega\right)  }}{\left\Vert b_{k-1,\mathbf{z}}\right\Vert
_{L^{2}\left(  \Omega\right)  }^{2}}.\label{SIIfinal}%
\end{align}
The conclusion of
\eqref{eqn:vfz_q_split}--\eqref{SIIfinal} reads
$\left(  \varphi_{\mathbf{z}},q\right)  _{L^{2}\left(
\Omega\right)  }=0$ and, since $\mathbf{z}\in \mathcal{S}\mathcal{C}_{\mathcal{T}}(\eta)$ was arbitrary,
$\mathcal{B}_0\subseteq M_{\eta,k-1}^{\operatorname{mod}%
}\left(  \mathcal{T}\right)  ^{\perp}$.
It remains to show that $\mathcal{B}_0$
is linearly independent. 
Given any $c_{\mathbf{z}}\in\mathbb{R}$ for $\mathbf{z}\in\mathcal S\mathcal{C}_{\mathcal{T}}(\eta)$ with 
$\sum_{\mathbf{z}\in\mathcal{SC}_{\mathcal{T}}\left(  \eta\right)
}c_{\mathbf{z}}\varphi_{\mathbf{z}}=0$, the definition of
$\varphi_{\mathbf{z}}$ in~\eqref{Eq:Def varphi} reveals
\begin{align*}
0&  =\sum_{\mathbf{z}\in\mathcal{SC}_{\mathcal{T}}\left(  \eta\right)
}\left(C_{\mathbf{z}}\phi_{\mathbf{z}}-c_{\mathbf{z}}\frac{\left(  b_{k-1,\mathbf{z}%
}\right)  _{\operatorname*{mvz}}}{\left\Vert b_{k-1,\mathbf{z}%
}\right\Vert _{L^{2}\left(  \Omega\right)  }^{2}}\right)%
\end{align*}
with $C_{\mathbf{z}}\coloneqq c_{\mathbf{z}}-\sum^{}_{\mathbf{z}\in\mathcal{S}\mathcal{C}_{\mathcal{T}}(\eta)}\left(
\overline{b_{k-1,\mathbf{z}}},b_{k-1,\mathbf{y}}\right)  _{L^{2}\left(
\Omega\right)  }/\left\Vert b_{k-1,\mathbf{z}%
}\right\Vert _{L^{2}\left(  \Omega\right)  }^{2}$. This leads to the condition%
\[
T:=\sum_{\mathbf{z}\in\mathcal{SC}_{\mathcal{T}}\left(  \eta\right)  }%
\frac{c_{\mathbf{z}}}{\left\Vert b_{k-1,\mathbf{z}}\right\Vert _{L^{2}\left(
\Omega\right)  }^{2}}\left(  b_{k-1,\mathbf{z}}\right)  _{\operatorname*{mvz}%
}=\sum_{\mathbf{y}\in\mathcal{SC}_{\mathcal{T}}\left(  \eta\right)
}C_{\mathbf{y}}\phi_{\mathbf{y}}\in M_{\eta,k-1}(\mathcal{T}).
\]
Hence, $T\in\operatorname{span}\{\left(  b_{k-1,\mathbf{z}}\right)  _{\operatorname*{mvz}%
}\ |\ \mathbf{z}\in \mathcal{S}\mathcal{C}_{\mathcal{T}}(\eta)\}\cap M_{\eta,k-1}(\mathcal{T})=\{0\}$ vanishes by
Lemma~\ref{PropOrthComplEta}.
Since the $\left(  b_{k-1,\mathbf{z}}\right)
_{\operatorname*{mvz}}$ are linear independent by Proposition~\ref{Propcritfct}, $T=0$ implies $c_{\mathbf{z}}=0$ for
all $\mathbf{z}\in\mathcal{S}\mathcal{C}_{\mathcal{T}}(\eta)$.
This verifies the linear independence of $\mathcal{B}_0$ and concludes the proof.
\end{proof}

The following lemma recalls the key estimates for the divergence control in \cite{Sauter_eta_wired}.

\begin{proposition}[{\cite{Sauter_eta_wired}}]
\label{PropDivControl}
(a)
Given $0\leq \eta\leq 1$, any $q\in M_{0,k-1}(\mathcal T)\cap M_{\eta,k-1}(\mathcal{T})^\perp$ satisfies
\begin{equation*}
\left\Vert q\right\Vert _{L^{2}\left(
\Omega\right)  }^2\leq\frac{12}{7}\binom{k+1}{2}^{-2}\sum
_{\mathbf{z}\in\mathcal{C}_{\mathcal{T}}\left(  \eta\right)  }h_{\mathbf{z}%
}^{2}\left(  A_{\mathcal{T},\mathbf{z}}\left(  q\right)  \right)  ^{2}. 
\end{equation*}
(b)
There exists $\eta_0>0$ exclusively depending on the shape-regularity of the mesh and the minimal outer angle
such that for any $0\leq \eta<\eta_0$ and any $\mathbf{z}\in\mathcal{C}_{\mathcal{T}}(\eta)$ we
have
\begin{align*}
		\left\vert A_{\mathcal{T}%
			,\mathbf{z}}\left(  \operatorname*{div}\mathbf{v}\right)  \right\vert \leq Ch_{\mathbf{z}}%
		^{-1}k^{2}\eta\left\Vert \nabla\mathbf{v}\right\Vert _{\mathbb{L}^{2}\left(
			\omega_{\mathbf{z}}\right)  }\quad\forall\mathbf{v}\in\mathbf{S}%
		_{k,0}\left(  \mathcal{T}\right)  .
\end{align*}
The constant $C>0$ depends only on the shape-regularity of the mesh.
%
\end{proposition}
\begin{proof}
The estimate in (a) is the conclusion of Step 2 of the proof of Lemma 4 in \cite{Sauter_eta_wired}, therein
stated for the divergence $\operatorname{div}\mathbf{v}_{\mathbf{S}}=q$ of the right-inverse $\mathbf{v}_{\mathbf{S}}=\Pi_kq\in \mathbf{S}_{k,0}(\mathcal{T})$ from
Proposition~\ref{Prop:InfSup stability mod SV}. The estimate in (b) is provided in \cite[Cor.~1 on p.~16]{Sauter_eta_wired} in terms of $\eta_0$ from
\cite[Lem.~2.13]{SAUTER202349} discussed in Remark~\ref{rem:eta0}; further details are omitted.
\end{proof}


\begin{proof}[Proof of Theorem \ref{Thm:Divergence estimate}]
The proof of Theorem~\ref{Thm:Divergence estimate} consists of three steps.

\textbf{Step 1: preparations.} 
Since $\operatorname{div}\mathbf{u}_{\mathbf{S}}\in M_{\eta,k-1}^{\operatorname{mod}}(\mathcal{T})^\perp$
by~\eqref{eqn:div_pw_orthogonality}, 
Lemma \ref{Prop:Orhtogonal complement Mmod}
provides
\begin{equation}
\operatorname{div}\mathbf{u}_{\mathbf{S}}=\sum_{\mathbf{z}\in\mathcal{C}%
_{\mathcal{T}}\left(  \eta\right)  \setminus\mathcal{SC}_{\mathcal{T}}\left(
\eta\right)  }c_{\mathbf{z}}b_{k-1,\mathbf{z}}+\sum_{\mathbf{y}\in
\mathcal{SC}_{\mathcal{T}}\left(  \eta\right)  }c_{\mathbf{y}}\varphi
_{\mathbf{y}}=q_{1}+q_{2}+q_{3}, \label{Eq:Decomposition div us}%
\end{equation}
for coefficients $c_{\mathbf{z}},c_{\mathbf{y}}, C_{\mathbf{y}}\in\mathbb{R}$ and
\begin{align*}
	q_{1}  :=\hspace*{-.7em}\sum_{\mathbf{z}\in\mathcal{C}_{\mathcal{T}}\left(  \eta\right)
\setminus\mathcal{SC}_{\mathcal{T}}\left(  \eta\right)  }\hspace*{-1em}c_{\mathbf{z}%
}b_{k-1,\mathbf{z}},\quad
q_{2}:=-\hspace*{-0.5em}\sum_{\mathbf{y}\in\mathcal{SC}_{\mathcal{T}}\left(  \eta\right)
}\frac{c_{\mathbf{y}}}{\left\Vert b_{k-1,\mathbf{y}}\right\Vert _{L^{2}\left(
\Omega\right)  }^{2}}\left(  b_{k-1,\mathbf{y}}\right)  _{\operatorname*{mvz}%
},\quad
q_{3}    :=\hspace*{-.3em}\sum_{\mathbf{y}\in\mathcal{SC}_{\mathcal{T}}\left(  \eta\right)
}C_{\mathbf{y}}
\phi_{\mathbf{y}}.
\end{align*}
More precisely, the coefficients $C_{\mathbf{y}}$ are given for all
$\mathbf{y}\in\mathcal{S}\mathcal{C}_{\mathcal{T}}(\eta)$ in terms of $c_{\mathbf{y}}$ as
\begin{align}\label{eqn:Cy_def}
	C_{\mathbf{y}} &:=c_{\mathbf{y}}-\sum_{\mathbf{a}\in\mathcal{SC}_{\mathcal{T}}\left(
\eta\right)  }c_{\mathbf{a}}\frac{\left(  \overline{b_{k-1,\mathbf{a}}%
},b_{k-1,\mathbf{y}}\right)  _{L^{2}\left(  \Omega\right)  }}{\left\Vert
b_{k-1,\mathbf{a}}\right\Vert _{L^{2}\left(  \Omega\right)  }^{2}}.
\end{align}
Since $\mathcal{SC}_{\mathcal{T}}\left(  \eta\right)  $
contains only Robinson vertices by assumption,
the functions 
$q_{1}$, $q_{2}$, $q_{3}$ are pairwise orthogonal by Lemma~\ref{PropOrthComplEta} and
$\phi_{\mathbf{z}}\in
M_{\eta,k-1}\left(  \mathcal{T}\right)  $ 
for all $\mathbf{z}\in\mathcal{C}_{\mathcal{T}}\left(  \eta\right)
\setminus\mathcal{SC}_{\mathcal{T}}\left(  \eta\right)  $.
Recall $\operatorname{div}%
\mathbf{u}_{\mathbf{S}}\in M_{0,k-1}\left(  \mathcal{T}\right)  $
from~\eqref{Eq:Def SV M} so that
$q_{1},q_{2}\in
M_{\eta,k-1}\left(  \mathcal{T}\right)  ^{\perp}$ and
$q_{3}\in M_{\eta,k-1}\left(  \mathcal{T}\right)
\subseteq M_{0,k-1}\left(  \mathcal{T}\right)  $ imply
\begin{align}\label{eqn:q1_q2_elem}
	\operatorname{div}\mathbf{u}_{\mathbf{S}}-q_{3}=q_1+q_2  \in M_{0,k-1}\left(
	\mathcal{T}\right) \cap
	M_{\eta,k-1}\left(  \mathcal{T}\right)  ^{\perp}.
\end{align}
\textbf{Step 2: control of }$q_{3}$\textbf{.} By assumption, \eqref{formfz}--\eqref{Eq:Jzcontestimate} hold for all $\mathbf{z}\in\mathcal{S}\mathcal{C}_{\mathcal{T}}(\eta)$.
In particular, \eqref{Eq:Jzcontestimate} bounds the operator norm 
$\left\Vert f_{\mathbf{z}}\right\Vert _{\mathbb{R}\leftarrow M_{\eta
,k-1}}$ of $f_{\mathbf{z}}$ and reveals with~\eqref{Eq:Riesz rep} that
\begin{align*}
	\left\Vert \phi_{\mathbf{z}%
		}\right\Vert _{L^{2}\left(  \Omega\right)  }=\left\Vert f_{\mathbf{z}%
			}\right\Vert _{\mathbb{R}\leftarrow M_{\eta,k-1}}\leq C_{f_{\mathbf{z}}%
			}/\left\Vert b_{k-1,\mathbf{z}}\right\Vert _{L^{2}\left(  \Omega\right)  }.
\end{align*}
This, triangle and Cauchy-Schwarz inequalities (in $\ell^2$) result for 
$\widetilde{C}_{f}:=\sqrt{\sum^{}_{\mathbf{y}\in\mathcal{S}\mathcal{C}_{\mathcal{T}}(\eta)}  C_{f_{\mathbf{y}}}^2}$ in
\begin{align}
\left\Vert q_{3}\right\Vert _{L^{2}\left(  \Omega\right)  }  &
\leq \sum^{}_{\mathbf{y}\in\mathcal{S}\mathcal{C}_{\mathcal{T}}(\eta)} \frac{C_{\mathbf{y}}}{\|b_{k-1,\mathbf{y}}\|_{L^2(\Omega)}}C_{f_{\mathbf{y}}}
\leq \widetilde{C}_f \sqrt{\sum^{}_{\mathbf{y}\in\mathcal{S}\mathcal{C}_{\mathcal{T}}(\eta)}
\frac{C_{\mathbf{y}}^2}{\|b_{k-1,\mathbf{y}}\|_{L^2(\Omega)}^2}}.
\end{align}
The orthogonality $\left\Vert \left(  b_{k-1,\mathbf{y}}\right)_{\operatorname*{mvz}}\right\Vert _{L ^{2}(\Omega)}^2+\left\Vert \overline{b_{k-1,\mathbf{z}}}\right\Vert _{L ^{2}\left(
\Omega\right)  }^2 = \left\Vert
{b_{k-1,\mathbf{z}}}\right\Vert _{L ^{2}\left(  \Omega\right)  }^2$ of the $L^{2}$-projection, the Cauchy-Schwarz
inequality, and
the definition of $C_{\mathbf{y}}$ in~\eqref{eqn:Cy_def} reveal 
\begin{align*}
	C_{\mathbf{y}}\leq c_{\mathbf{y}}\frac{\|\left(  b_{k-1,\mathbf{y}}\right)  _{\operatorname*{mvz}%
}\|_{L^2(\Omega)}}{\left\Vert
b_{k-1,\mathbf{y}}\right\Vert _{L^{2}\left(  \Omega\right)  }^{}}+\sum_{\mathbf{a}\in\mathcal{SC}_{\mathcal{T}}\left(
	\eta\right)\setminus\{\mathbf{y}\}  }c_{\mathbf{a}}\frac{\|{b_{k-1,\mathbf{y}}%
}\|_{L^2(\Omega)}}{\left\Vert
b_{k-1,\mathbf{a}}\right\Vert _{L^{}\left(  \Omega\right)  }^{2}}
\leq \sum_{\mathbf{a}\in\mathcal{SC}_{\mathcal{T}}\left(
	\eta\right)}c_{\mathbf{a}}\frac{\left\Vert
b_{k-1,\mathbf{y}}\right\Vert _{L^{2}\left(  \Omega\right)  }}{\left\Vert
b_{k-1,\mathbf{a}}\right\Vert _{L^{2}\left(  \Omega\right)  }}.
\end{align*}
The previous two estimates and a Cauchy-Schwarz inequality imply for
$C_{0}=\widetilde{C}_f^2\operatorname{card}\mathcal{SC}_{\mathcal{T}}\left(  \eta\right)^2$
that
\begin{align}\label{eqn:q3_est}
\left\Vert q_{3}\right\Vert _{L^{2}\left(  \Omega\right)  }^{2}\leq 
\widetilde{C}_f^2 \sum^{}_{\mathbf{y}\in\mathcal{S}\mathcal{C}_{\mathcal{T}}(\eta)}
\left(\sum^{}_{\mathbf{a}\in\mathcal{S}\mathcal{C}_{\mathcal{T}}(\eta)}
\frac{c_{\mathbf{a}}}{\|b_{k-1,\mathbf{a}}\|_{L^2(\Omega)}}\right)^2
\leq
C_{0}%
^{2}\sum_{\mathbf{y}\in\mathcal{SC}_{\mathcal{T}}\left(  \eta\right)  }%
\frac{c_{\mathbf{y}}^{2}}{\left\Vert b_{k-1,\mathbf{y}}\right\Vert
_{L^{2}\left(  \Omega\right)  }^{2}}.
\end{align}
All $\eta$-super critical vertices are of Robinson
type by assumption and
$\left\{  \left.  b_{k-1,\mathbf{y}}\;\right\vert
\;\mathbf{y}\in\mathcal{SC}_{\mathcal{T}}\left(  \eta\right)  \right\}  $
is pairwise orthogonal by Remark~\ref{Rem:Avoiding techinacl difficulties}.(ii). Hence, the sum in the right-hand side of~\eqref{eqn:q3_est} is the $L^2$ norm of 
\[
\tilde{q}_{2}:=-\sum_{\mathbf{y}\in\mathcal{SC}_{\mathcal{T}}\left(
\eta\right)  }\frac{c_{\mathbf{y}}}{\left\Vert b_{k-1,\mathbf{y}}\right\Vert
_{L^{2}\left(  \Omega\right)  }^2}b_{k-1,\mathbf{y}}.
\]
A comparison with the definition of $q_2$ verifies $q_{2}=\left(  \tilde{q}_{2}\right)
_{\operatorname*{mvz}}$ and \eqref{stabest} in Proposition \ref{Propcritfct}
shows%
\[
C_{1}^{-2}\|q_3\|_{L^2(\Omega)}^2\leq \left\Vert \tilde{q}_{2}\right\Vert _{L^{2}\left(  \Omega\right)  }^{2}%
\leq\frac{16}{7}\sum_{K\in\mathcal{T}}\inf_{\alpha\in\mathbb{R}}\left\Vert
\tilde{q}_{2}-\alpha\right\Vert _{L^{2}\left(  K\right)  }^{2}\leq\frac{16}%
{7}\left\Vert q_{2}\right\Vert _{L^{2}\left(  \Omega\right)  }^{2}.
\]
The conclusion of the previous estimates is the existence of $C>0$ exclusively depending on the cardinality of
$\mathcal{SC}_{\mathcal{T}}\left(  \eta\right)  $ and the constants
$C_{f_{\mathbf{y}}}$ such that
\begin{equation}
\left\Vert q_{3}\right\Vert _{L^{2}\left(  \Omega\right)  }\leq C\left\Vert
q_{2}\right\Vert _{L^{2}\left(  \Omega\right)  }.
\label{Eq:Control q_3 by q_2}%
\end{equation}
\textbf{Step 3: Conclusion.} The pairwise orthogonality of $q_1,q_2,q_3$ with $\operatorname*{div}\mathbf{u}_{\mathbf{S}}=q_1+q_2+q_3$ from Step 1
and
\eqref{Eq:Control q_3 by q_2} verify
\begin{align}
\left\Vert \operatorname{div}\mathbf{u}_{\mathbf{S}}\right\Vert _{L^{2}\left(
\Omega\right)  }^{2}  &  =\left\Vert q_{1}\right\Vert _{L^{2}\left(
\Omega\right)  }^{2}+\left\Vert q_{2}\right\Vert _{L^{2}\left(  \Omega\right)
}^{2}+\left\Vert q_{3}\right\Vert _{L^{2}\left(  \Omega\right)  }^{2}%
\leq\left\Vert q_{1}\right\Vert _{L^{2}\left(  \Omega\right)  }^{2}+\left(
1+C\right)  \left\Vert q_{2}\right\Vert _{L^{2}\left(  \Omega\right)  }%
^{2}\nonumber\\
&  \leq\left(  1+C\right)  \left\Vert q_{1}+q_{2}\right\Vert _{L^{2}\left(
\Omega\right)  }^{2}. \label{Eq:Estimate div us by q1 q2}%
\end{align}
Proposition \ref{PropDivControl}.(a) with~\eqref{eqn:q1_q2_elem} and $A_{\mathcal{T},\mathbf{z}}(q_3)=0$ from $q_3\in M_{\eta,k-1}(\mathcal{T})$
provide
\begin{align}\label{eqn:q1_q2_control}
\left\Vert q_{1}+q_{2}\right\Vert _{L^{2}\left(  \Omega\right)  }^{2}
\leq
C\sum_{\mathbf{z}\in\mathcal{C}_{\mathcal{T}}\left(  \eta\right)
}h_{\mathbf{z}}^{2}\left[  A_{\mathcal{T},\mathbf{z}}\left(  q_{1}%
+q_{2}\right)  \right]  ^{2}
\leq
C\sum_{\mathbf{z}\in\mathcal{C}_{\mathcal{T}}\left(  \eta\right)
}h_{\mathbf{z}}^{2}\left[  A_{\mathcal{T},\mathbf{z}}\left(  \operatorname{div}\mathbf{u}_{\mathbf{S}}\right)  \right]  ^{2}.%
\end{align}
By definition of $\mathbf{S}_{\eta,k,0}$ in~\eqref{DefSetak0T},
$A_{\mathcal{T},\mathbf{y}}(\operatorname*{div}\mathbf{w}_{\mathbf{S}}) =0$ vanishes for all
$\mathbf{w}_{\mathbf{S}}\in \mathbf{S}_{\eta,k,0}$ and $\mathbf{y}\in
\mathcal{C}_{\mathcal{T}}(\eta)$.
Proposition~\ref{PropDivControl}.(b) provides $0<\eta_0$ such that this and the combination of~\eqref{Eq:Estimate div us by q1 q2}--\eqref{eqn:q1_q2_control} 
reveal for $0\leq \eta< \eta_0$ that
\[
\left\Vert \operatorname{div}\mathbf{u}_{\mathbf{S}}\right\Vert _{L^{2}\left(
\Omega\right)  }^{2}\leq C\sum_{\mathbf{z}\in\mathcal{C}_{\mathcal{T}}\left(
\eta\right)  }h_{\mathbf{z}}^{2}A_{\mathcal{T},\mathbf{z}}\left(
\operatorname{div}\left(  \mathbf{u}_{\mathbf{S}}-\mathbf{w}_{\mathbf{S}%
}\right)  \right)  ^{2}\leq C\eta^{2}\left\Vert \nabla\left(  \mathbf{u}%
_{\mathbf{S}}-\mathbf{w}_{\mathbf{S}}\right)  \right\Vert _{\mathbb{L}%
^{2}\left(  \Omega(\eta)\right)  }^{2}
\]
using the finite overlay of the vertex patches $\{\omega_{\mathbf{z}}\ |\ \mathbf{z}\in
\mathcal{C}_{\mathcal{T}}(\eta)\}$ in the neighbourhood $\Omega(\eta)$ from~\eqref{eqn:Omega_star_def}. This concludes the proof with $\eta_0$ from Proposition~\ref{PropDivControl}.(b).
\end{proof}

With this we obtain an analogous estimates as in \cite[Thm. 3]{Sauter_eta_wired}.

\begin{proposition}
\label{Prop:Almost pressure robust} Under the Assumptions of Theorem~\ref{Thm:Divergence estimate} with the solution $\left(  \mathbf{u},p\right)  $ to
\eqref{varproblemstokes} for $\mathbf{F}\in\mathbf{H}^{-1}\left(
	\Omega\right)  $, the discrete solution $\left(  \mathbf{u}_{\mathbf{S}}%
,p_{M}^{\ast}\right)  \in\left(  \mathbf{S}_{k,0}\left(  \mathcal{T}\right)
,M_{\eta,k-1}^{\operatorname{mod}}\left(  \mathcal{T}\right)  \right)  $ to
\eqref{discrStokes} satisfies%
\begin{subequations}
\label{Propests}
\end{subequations}%
\begin{align}
\left\Vert \nabla\left(  \mathbf{u}-\mathbf{u}_{\mathbf{S}}\right)
\right\Vert _{\mathbb{L}^{2}\left(  \Omega\right)  }  &  \leq C\left(
\inf_{\substack{\mathbf{v}_{\mathbf{S}}\in\mathbf{S}_{k,0}\left(
\mathcal{T}\right)  \\\operatorname{div}\mathbf{v}_{\mathbf{S}}=0}}\left\Vert
\nabla\left(  \mathbf{u}-\mathbf{v}_{\mathbf{S}}\right)  \right\Vert
_{\mathbb{L}^{2}\left(  \Omega\right)  }+\eta\inf_{q_{M}\in M_{\eta
,k-1}^{\operatorname*{mod}}\left(  \mathcal{T}\right)  }\left\Vert
p-q_{M}\right\Vert _{L^{2}\left(  \Omega\right)  }\right)  ,\tag{%
\ref{Propests}%
a}\label{Propestsa1}\\
\left\Vert \operatorname{div}\mathbf{u}_{\mathbf{S}}\right\Vert _{L^{2}\left(
\Omega\right)  }  &  \leq C\left(  \eta\inf_{\substack{\mathbf{v}_{\mathbf{S}%
}\in\mathbf{S}_{k,0}\left(  \mathcal{T}\right)  \\\operatorname{div}%
\mathbf{v}_{\mathbf{S}}=0}}\left\Vert \nabla\left(  \mathbf{u}-\mathbf{v}%
_{\mathbf{S}}\right)  \right\Vert _{\mathbb{L}^{2}\left(  \Omega\right)
}+\eta^{2}\inf_{q_{M}\in M_{\eta,k-1}^{\operatorname*{mod}}\left(
\mathcal{T}\right)  }\left\Vert p-q_{M}\right\Vert _{L^{2}\left(
\Omega\right)  }\right)  . \tag{%
\ref{Propests}%
b}\label{Propestsb1}%
\end{align}

\end{proposition}
\begin{proof}
We start with the first estimate \eqref{Propestsa1} and follow the arguments
of the proof of \cite[Thm. 3]{Sauter_eta_wired}. Since the discrete velocity
$\mathbf{u}_{\mathbf{S}}$ satisfies the second equation in (\ref{discrStokes})
we conclude \eqref{eqn:div_pw_orthogonality}. Next,
consider some $\mathbf{v}_{\mathbf{S}}\in\mathbf{S}_{k,0}\left(
\mathcal{T}\right)  $ with $\operatorname*{div}\mathbf{v}_{\mathbf{S}}=0$ and observe that $\mathbf{v}_{\mathbf{S}}\in\mathbf{S}_{\eta,k,0}\left(
\mathcal{T}\right)  $, such that the first equations in
\eqref{varproblemstokes} and \eqref{discrStokes} for the test function
$\mathbf{e}_{\mathbf{S}}=\mathbf{u}_{\mathbf{S}}-\mathbf{v}_{\mathbf{S}}$
verifies
\[
a\left(  \mathbf{u}-\mathbf{u}_{\mathbf{S}},\mathbf{e}_{\mathbf{S}}\right)
=b\left(  \mathbf{e}_{\mathbf{S}},p-p_{M}\right)  .
\]
Since $\operatorname{div}\mathbf{e}_{\mathbf{S}}=\operatorname*{div}%
\mathbf{u}_{\mathbf{S}}$ is orthogonal to $M_{\eta,k-1}^{\operatorname*{mod}%
}\left(  \mathcal{T}\right)  $ we get%
\[
a\left(  \mathbf{u}-\mathbf{u}_{\mathbf{S}},\mathbf{e}_{\mathbf{S}}\right)
=b\left(  \mathbf{u}_{\mathbf{S}},p-q_{M}\right)  \qquad\forall q_{M}\in
M_{\eta,k-1}^{\operatorname*{mod}}\left(  \mathcal{T}\right)  .
\]
The combination of these relations leads for any $q_{M}\in
M_{\eta,k-1}^{\operatorname*{mod}}\left(  \mathcal{T}\right)$ to%
\begin{align*}
\left\Vert \nabla\left(  \mathbf{u}_{\mathbf{S}}-\mathbf{v}_{\mathbf{S}%
}\right)  \right\Vert _{\mathbb{L}^{2}\left(  \Omega\right)  }^{2} &
=a\left(  \mathbf{u}_{\mathbf{S}}-\mathbf{v}_{\mathbf{S}},\mathbf{e}%
_{\mathbf{S}}\right)  =a\left(  \mathbf{u}-\mathbf{v}_{\mathbf{S}}%
,\mathbf{e}_{\mathbf{S}}\right)  -a\left(  \mathbf{u}-\mathbf{u}_{\mathbf{S}%
},\mathbf{e}_{\mathbf{S}}\right)  \\
&  \leq\left\Vert \nabla\left(  \mathbf{u}-\mathbf{v}_{\mathbf{S}}\right)
\right\Vert _{\mathbb{L}^{2}\left(  \Omega\right)  }\left\Vert \nabla
\mathbf{e}_{\mathbf{S}}\right\Vert _{\mathbb{L}^{2}\left(  \Omega\right)
}+\left\Vert \operatorname*{div}\mathbf{u}_{\mathbf{S}}\right\Vert
_{L^{2}\left(  \Omega\right)  }\left\Vert p-q_{M}\right\Vert _{L^{2}\left(
\Omega\right)  }.
\end{align*}
From (\ref{divest}) we obtain for the choice $\mathbf{w}_{\mathbf{S}%
}=\mathbf{v}_{\mathbf{S}}$, that
\begin{align} \label{Eq:Est div us via vs}
\left\Vert \operatorname{div}\mathbf{u}_{\mathbf{S}}\right\Vert _{L^{2}\left(
\Omega\right)  }\leq C_{\operatorname*{div}}\eta\left\Vert \nabla\left(
\mathbf{u}_{\mathbf{S}}-\mathbf{v}_{\mathbf{S}}\right)  \right\Vert
_{L^{2}\left(  \Omega\right)  }%
\end{align}
and, in turn,%
\[
\left\Vert \nabla\left(  \mathbf{u}_{\mathbf{S}}-\mathbf{v}_{\mathbf{S}%
}\right)  \right\Vert _{\mathbb{L}^{2}\left(  \Omega\right)  }\leq\left\Vert
\nabla\left(  \mathbf{u}-\mathbf{v}_{\mathbf{S}}\right)  \right\Vert
_{\mathbb{L}^{2}\left(  \Omega\right)  }+C_{\operatorname*{div}}\eta\left\Vert
p-q_{M}\right\Vert _{L^{2}\left(  \Omega\right)  }.
\]
Thus, the velocity error can be estimated by%
\begin{align*}
\left\Vert \nabla\left(  \mathbf{u}-\mathbf{u}_{\mathbf{S}}\right)
\right\Vert _{\mathbb{L}^{2}\left(  \Omega\right)  } &  \leq\left\Vert
\nabla\left(  \mathbf{u}-\mathbf{v}_{\mathbf{S}}\right)  \right\Vert
_{\mathbb{L}^{2}\left(  \Omega\right)  }+\left\Vert \nabla\left(
\mathbf{u}_{\mathbf{S}}-\mathbf{v}_{\mathbf{S}}\right)  \right\Vert
_{\mathbb{L}^{2}\left(  \Omega\right)  }\nonumber\\
&  \leq2\left\Vert \nabla\left(  \mathbf{u}-\mathbf{v}_{\mathbf{S}}\right)
\right\Vert _{\mathbb{L}^{2}\left(  \Omega\right)  }+C_{\operatorname*{div}%
}\eta\inf_{q_{M}\in M_{\eta,k-1}^{\operatorname*{mod}}\left(  \mathcal{T}%
\right)  }\left\Vert p-q_{M}\right\Vert _{L^{2}\left(  \Omega\right)
}.
\end{align*}
Since $\mathbf{v}_{\mathbf{S}}$ was chosen arbitrarily, we obtain \eqref{Propestsa1} for $C = 2 + C_{\operatorname{div}}$. For \eqref{Propestsb1}, we take \eqref{Eq:Est div us via vs} and obtain through a triangle inequality that
\begin{align*}
\left\Vert \operatorname{div}\mathbf{u}_{\mathbf{S}}\right\Vert _{L^{2}\left(
\Omega\right)  }
&\leq 
C_{\operatorname*{div}} \eta
\left( \left\Vert \nabla \left( \mathbf{u} - \mathbf{u}_{\mathbf{S}} \right) \right\Vert_{L ^{2}\left(  \Omega\right)  } 
+
\left\Vert \nabla \mathbf{u} - \mathbf{v}_{\mathbf{S}} \right\Vert_{L ^{2}\left(  \Omega\right)  }\right).
\end{align*}
This and \eqref{Propestsa1} reveals \eqref{Propestsb1} concluding the proof.
\end{proof}

Two options for obtaining convergence rates for the discrete solution from the
quasi-optimal error estimates are sketched in the following
remark.

\begin{remark}
	\label{Rem:Explicite rates}Let $\left(  \mathbf{u}_{\mathbf{S}},p_{M}^{\ast}\right)
\in\left(  \mathbf{S}_{k,0}\left(  \mathcal{T}\right)  ,M_{\eta,k-1}%
^{\operatorname{mod}}\left(  \mathcal{T}\right)  \right)  $ denote the
discrete solution of the Stokes problem.

\textbf{1. }The combination of Corollary \ref{Cor:Best approx pwS mod} with
Theorem \ref{Thm:Recoverd approx Meta mod} directly results in $h$-explicit
convergence rates for the infima in the quasi-optimal error estimates
(\ref{Corests}). A minor drawback is that in this way the pressure infima are
not multiplied by the small pre-factor $\eta$ in (\ref{Corests}) in contrast
to (\ref{Propests}), i.e., the resulting estimate is less
\emph{pressure-robust.}

\textbf{2. }A direct estimate of the best-approximation of the velocity in
(\ref{Propests}) can be obtained by selecting for $\mathbf{v}_{\mathbf{S}}$
the original (solenoidal) Scott-Vogelius velocity field $\mathbf{v}%
_{\mathbf{S}}^{\operatorname{SV}}$. However, the stability bound of
$\mathbf{v}_{\mathbf{S}}^{\operatorname{SV}}$ is reciprocally related to
discrete inf-sup constant of the original Scott-Vogelius
element and may become large if the mesh has nearly-singular vertices.

A direct construction of some \emph{solenoidal} $\mathbf{v}_{\mathbf{S}}%
\in\mathbf{S}_{k,0}\left(  \mathcal{T}\right)  $ which leads to optimal,
mesh-robust convergence rates for the first infima in (\ref{Propests}) is, to
the best of our knowledge, an open question.
\end{remark}

\appendix

\section{Approximation by polynomial extension}

The proofs for the approximation property of our modified pressure spaces rely
on some estimates for the analytic extension of polynomials (extension by
\textquotedblleft itself\textquotedblright) proved in this
appendix. For a triangle $K$ with barycenter $\mathbf{M}_{K}$ we introduce a
neighborhood in terms of some given $\alpha\in[0,2\pi[$ and a scaling parameter $\lambda\geq0$ below. Consider the infinite line
\begin{equation}
L_{\alpha}:=\left\{  \mathbf{M}_{K}+s\binom{\cos\alpha}{\sin\alpha}%
:s\in\mathbb{R}\right\}  \text{.} \label{defLalpha}%
\end{equation}
and let $s_{+}>0$ be such that $\mathbf{y}_{\alpha}:=\mathbf{M}_{K}+s_{+}\binom{\cos\alpha}{\sin\alpha
}\in L_{\alpha}\cap\partial K$. 
Denote the midpoint of the intersection $L_{\alpha}\cap K$ by
$\mathbf{M}_{\alpha}\in K$. Given $\lambda\geq0$, the neighbourhood reads
\begin{equation}
K_{\lambda}:=\left\{  \mathbf{M}_{\alpha}+s\left(  \mathbf{y}_{\alpha
}-\mathbf{M}_{\alpha}\right)  \mid\left\{
\begin{array}
[c]{l}%
0\leq s\leq1+\lambda\\
0\leq\alpha<2\pi
\end{array}
\right.  \quad\right\}  . \label{defklambda}%
\end{equation}
Clearly $K=K_{0}\subseteq K_{\lambda}$ for $\lambda\geq0$.

\begin{lemma}
\label{Lemy}There exists $c>0$ depending solely on the shape regularity of $K$
such that for any $\lambda\geq0:$%
\begin{equation}
\left\{  \mathbf{y}\in\mathbb{R}^{2}\backslash K\mid c\frac
{\operatorname{dist}\left(  \mathbf{y},K\right)  }{h_{K}}\leq\lambda\right\}
\subset K_{\lambda}\quad\text{and\quad}\operatorname*{diam}K_{\lambda}%
\leq\left(  1+2\lambda\right)  h_{K}. \label{Krincl}%
\end{equation}

\end{lemma}%

\begin{proof}
Since $K_{\lambda}$ is compact, there exist
$\mathbf{y},\mathbf{z}\in\partial K_{\lambda}$ such that
\begin{equation}
\operatorname{diam}K_{\lambda}=\left\Vert \mathbf{y}-\mathbf{z}\right\Vert
.\label{Eq:Classification of diameter}%
\end{equation}
From \eqref{defklambda} it follows that there exist $0\leq\alpha,\beta<2\pi$
such that
\[
\mathbf{y}=\mathbf{M}_{\alpha}+\left(  1+\lambda\right)  \left(
\mathbf{y}_{\alpha}-\mathbf{M}_{\alpha}\right)  \quad\text{and}\quad
\mathbf{z}=\mathbf{M}_{\beta}+\left(  1+\lambda\right)  \left(  \mathbf{z}%
_{\beta}-\mathbf{M}_{\beta}\right)  .
\]
Rearranging both terms yields,
\begin{equation}
\mathbf{y}=\mathbf{y}_{\alpha}+\lambda\left(  \mathbf{y}_{\alpha}%
-\mathbf{M}_{\alpha}\right)  \quad\text{and}\quad\mathbf{z}=\mathbf{z}_{\beta
}+\lambda\left(  \mathbf{z}_{\beta}-\mathbf{M}_{\beta}\right)
.\label{Eq:reaaranged outer terms}%
\end{equation}
Two triangle inequalities and 
\eqref{Eq:Classification of diameter}--\eqref{Eq:reaaranged outer terms} result in
\begin{align*}
\operatorname{diam}K_{\lambda} &  \leq\left\Vert \mathbf{y}-\mathbf{y}%
_{\alpha}\right\Vert +\left\Vert \mathbf{y_{\alpha}-\mathbf{z}_{\beta}%
}\right\Vert +\left\Vert \mathbf{z}-\mathbf{z}_{\beta}\right\Vert \leq
h_{K}+\lambda\left(  \left\Vert \mathbf{y}_{\alpha}-\mathbf{M}_{\alpha
}\right\Vert +\left\Vert \mathbf{z}_{\beta}-\mathbf{M}_{\beta}\right\Vert
\right)  \\
&  \leq h_{K}+2\lambda\sup_{0\leq\gamma<2\pi}\left\Vert \mathbf{y}_{\gamma
}-\mathbf{M}_{\gamma}\right\Vert \leq\left(  1+2\lambda\right)  h_{K}.
\end{align*}
To prove the first inclusion in (\ref{Krincl}), let $\mathbf{y}\in
\mathbb{R}^{2}\backslash K$ and $\alpha\in\left[  0,2\pi\right[  $, $s\geq0$
such that%
\begin{equation}
\mathbf{y}=\mathbf{M}_{\alpha}+\left(  1+s\right)  \left(  \mathbf{y}_{\alpha
}-\mathbf{M}_{\alpha}\right)  .\label{defyma}%
\end{equation}
Rearranging the terms as in \eqref{Eq:reaaranged outer terms} implies that
\[
\left\Vert \mathbf{y}-\mathbf{y}_{\alpha}\right\Vert =s\left\Vert
\mathbf{y}_{\alpha}-\mathbf{M}_{\alpha}\right\Vert .
\]
Elementary trigonometry provides a constant $0<c$ depending only on the shape regularity of $K$
such that $\left\Vert \mathbf{y}-\mathbf{y}_{\alpha}\right\Vert \leq
\sqrt c\operatorname{dist}\left(  \mathbf{y},K\right)  $ and $\left\Vert
\mathbf{y}_{\alpha}-\mathbf{M}_{\alpha}\right\Vert \geq h_{K}/\sqrt c$ so that%
\[
s=\frac{\left\Vert \mathbf{y}-\mathbf{y}_{\alpha}\right\Vert }{\left\Vert
\mathbf{y}_{\alpha}-\mathbf{M}_{\alpha}\right\Vert }\leq c%
\frac{\operatorname{dist}\left(  \mathbf{y},K\right)  }{h_{K}}.
\]
If the right-hand side is bounded by $\lambda$, then $s\leq \lambda$
follows and~\eqref{defyma} implies $\mathbf{y}\in K_{\lambda}$. This concludes the proof.
\end{proof}

\begin{lemma}
	\label{Lemextension} Given a triangle $K$, there exists $c\geq 1$ depending only on the shape regularity of $K$ such
	that any
polynomial $p\in\mathbb{P}_{k}\left(  \mathbb{R}^{2}\right)  $ satisfies%
\begin{align} \label{Eq:Estimate extended polynomial values via Linfty norm}
\left\vert p\left(  \mathbf{y}\right)  \right\vert \leq\left\vert T_{k}\left(
1+\lambda\right)  \right\vert \left\Vert p\right\Vert _{L^{\infty}\left(
K\right)  }\quad\forall\mathbf{y}\in\mathbb{R}^{2}\backslash K\quad\text{for
}\lambda:=c\frac{\operatorname{dist}\left(  \mathbf{y},K\right)  }{h_{K}}.%
\end{align} 
The Chebyshev
polynomial $T_{k}$ outside the interval $\left[  -1,1\right]  $ admits
\begin{equation}
\left\vert T_{k}\left(  1+\lambda\right)  \right\vert \leq\gamma_{k}\left(
\lambda\right)  :=\min\left\{  \frac{2^{k}+1}{2}\left(  1+\lambda\right)
^{k},\operatorname*{e}\nolimits^{\lambda k^{2}/2}\right\}  . \label{defgammak}%
\end{equation}

\end{lemma}%

\begin{proof}
We start with a one-dimensional consideration. It is well known (see, e.g.,
\cite[Prop.~2.3, 2.4]{sachdeva2014faster}) that any polynomial $f\in
\mathbb{P}_{k}\left(  \mathbb{R}\right)  $ satisfies%
\[
\left\vert f\left(  y\right)  \right\vert \leq\left\Vert f\right\Vert
_{L^{\infty}\left(  \left[  -1,1\right]  \right)  }\max\left\{  1,\left\vert
T_{k}\left(  y\right)  \right\vert \right\}  \quad\forall y\in\mathbb{R}.
\]
An affine transformation to the interval $I=\left[  m-\varepsilon
,m+\varepsilon\right]  $ for $m\in\mathbb{R}$ and $\varepsilon>0$ provides
\begin{equation}
\left\vert f\left(  y\right)  \right\vert \leq\left\Vert f\right\Vert
_{L^{\infty}\left(  I\right)  }\max\left\{  1,\left\vert T_{k}\left(
\frac{y-m}{\varepsilon}\right)  \right\vert \right\}  \quad\forall
y\in\mathbb{R}. \label{fymm}%
\end{equation}
Next, consider $p\in\mathbb{P}_{k}\left(  \mathbb{R}^{2}\right)  $. 
Any $\mathbf{y}\in K_{\lambda}\backslash K$ fulfills $\mathbf{y}=\mathbf{M}_{\alpha
}+\left(  1+s\right)  \left(  \mathbf{y}_{\alpha}-\mathbf{M}_{\alpha}\right)
$ for some $\alpha\in\left[  0,2\pi\right[  $, $0\leq s\leq\lambda$ with
$\mathbf{y}_{\alpha}\in\partial K$ and $\mathbf{M}_{\alpha}\in K$. We choose
a coordinate system such that $\mathbf{M}_{\alpha}$ is the origin and
$\mathbf{y}=\left(  y_{1},0\right)  ^{T}$ for some $y_{1}>0$ and set
$\rho:=\left\Vert \mathbf{y}_{\alpha}-\mathbf{M}_{\alpha}\right\Vert $. 
Estimate (\ref{fymm}) applies to $f=p|_{y_2=0}\in \mathbb P_k(\mathbb R)$ for the interval $\left[
-\rho,\rho\right]  \times\left\{  0\right\}  \subset K$ with $y-m=\left\Vert \mathbf{y}-\mathbf{M}_{\alpha}\right\Vert
=\left(  1+s\right)  \rho$. Since $s\geq0$, it follows that
\[
\left\vert p\left(  \mathbf{y}\right)  \right\vert \leq\left\Vert p\right\Vert
_{L^{\infty}\left(  K\right)  }T_{k}\left(  1+s\right)  \leq\left\Vert
p\right\Vert _{L^{\infty}\left(  K\right)  }T_{k}\left(  1+\lambda\right)  .
\]
Finally, we prove (\ref{defgammak}). For $y\in\mathbb{R}$ with $\left\vert
y\right\vert \geq1$, the Chebyshev polynomials have the representation (see,
e.g., \cite[Prop. 2.5]{sachdeva2014faster})%
\[
T_{k}\left(  y\right)  =\frac{\left(  y+\sqrt{y^{2}-1}\right)  ^{k}+\left(
y-\sqrt{y^{2}-1}\right)  ^{k}}{2}.
\]
From this $\left\vert T_{k}\left(  y\right)  \right\vert
\leq \left(2^k +1\right)/2\left\vert y\right\vert ^{k}$ follows directly for $\left\vert y\right\vert \geq1$.
Hence,%
\begin{equation}
\left\vert p\left(  \mathbf{y}\right)  \right\vert \leq\left\Vert p\right\Vert
_{L^{\infty}\left(  K\right)  }\frac{2^{k}+1}{2}\left(  1+\lambda\right)
^{k}\quad\forall\mathbf{y}\in K_{\lambda}. \label{estf}%
\end{equation}
If $\mathbf{y}$ is close to $K$, this estimate can be improved. From \cite[Lem.
A.1]{grote2021stabilized} we conclude that%
\[
\left\vert T_{k}\left(  y\right)  \right\vert \leq\operatorname*{e}%
\nolimits^{\left(  \left\vert y\right\vert -1\right)  k^{2}/2}\quad\text{for
}\left\vert y\right\vert \geq1.
\]
Consequently, we have
\[
\left\vert p\left(  \mathbf{y}\right)  \right\vert \leq\left\Vert p\right\Vert
_{L^{\infty}\left(  K\right)  }\operatorname*{e}\nolimits^{\lambda k^{2}%
/2}\quad\forall\mathbf{y}\in K_{\lambda}.
\]
From Lemma \ref{Lemy}, the estimate \eqref{Eq:Estimate extended polynomial values via Linfty norm} for the relevant values of $\mathbf{y}$ follows.
\end{proof}

The next lemma is concerned with a neighborhood property for polynomial approximations.

\begin{lemma}
\label{LemNeighborhood}Let $p\in H^{s-1}\left(  \mathbb{R}^{2}\right)  $ for
some $s>1$. Given two triangles $K, K^*$, let $\delta\geq 0$ satisfy%
\[
K\subset K^{\ast}\subset\left\{  \mathbf{y}\in\mathbb{R}^{2}\mid
\operatorname{dist}\left(  \mathbf{y},K\right)  \leq\delta h_{K}\right\}.
\]
Then, for any polynomial $p_{k}\in\mathbb{P}_{k}\left(  \mathbb{R}^{2}\right)
$, it holds%
\[
\left\Vert p-p_{k}\right\Vert _{L^{2}\left(  K^{\ast}\right)  }\leq
T_{k}\left(  1+c\delta\right)  \left(  2\frac{\left(  \left(  1+2\delta
\right)  h_{K}\right)  ^{\min\left\{  k+1,s-1\right\}  }}{\left(  k+1\right)
^{s-1}}\left\Vert p\right\Vert _{H^{s-1}\left(  K^{\ast}\right)  }+\left\Vert
p-p_{k}\right\Vert _{L^{2}\left(  K\right)  }\right)  ,
\]
where $c$ and $C$ only depend on the shape regularity of $K$ and $K^{\ast}$.
The Chebyshev polynomial $T_{k}$ can be estimated by $\gamma_{k}\left(
c\delta\right)  $ with $\gamma_{k}$ as in (\ref{defgammak}).
\end{lemma}

\begin{proof}
We employ the operator $W_{k}:H^{s-1}\left(  K^{\ast}\right)  \rightarrow
\mathbb{P}_{k}\left(  K^{\ast}\right)  $ from \cite[Cor. 3.2]%
{Ainsworth_parker_II} based on the operator in \cite[Lem. 3.1]%
{Suri_hp_interpol} to obtain%
\begin{equation}
\left\Vert p-W_{k}p\right\Vert _{L^{2}\left(  K^{\ast}\right)  }\leq
C\frac{\left(  \left(  1+2\delta\right)  h_{K}\right)  ^{\min\left\{
k+1,s-1\right\}  }}{\left(  k+1\right)  ^{s-1}}\left\Vert p\right\Vert
_{H^{s-1}\left(  K^{\ast}\right)  }. \label{extapprox}%
\end{equation}
A triangle inequality leads to%
\[
\left\Vert p-p_{k}\right\Vert _{L^{2}\left(  K^{\ast}\right)  }\leq\left\Vert
p-W_{k}p\right\Vert _{L^{2}\left(  K^{\ast}\right)  }+\left\Vert W
_{k}p-p_{k}\right\Vert _{L^{2}\left(  K^{\ast}\right)  }.
\]
For any $p_k\in \mathbb P_k(K^*)$, Lemma \ref{Lemextension} with $\lambda
=c\delta$ provides
\begin{align*}
\left\Vert W_{k}p-p_{k}\right\Vert _{L^{2}\left(  K^{\ast}\right)  }  &
=T_{k}\left(  1+\lambda\right)  \left\Vert \left(  W_{k}p-p_{k}\right)
\right\Vert _{L^{2}\left(  K\right)  }\\
&  \leq T_{k}\left(  1+\lambda\right)  \left(  \left\Vert p-W_{k}%
p\right\Vert _{L^{2}\left(  K^{\ast}\right)  }+\left\Vert p-p_{k}\right\Vert
_{L^{2}\left(  K\right)  }\right)  ,
\end{align*}
with a triangle inequality in the last step.
This and a triangle inequality lead to%
\[
\left\Vert p-p_{k}\right\Vert _{L^{2}\left(  K^{\ast}\right)  }\leq\left(
1+T_{k}\left(  1+\lambda\right)  \right)  \left\Vert p-W_{k}p\right\Vert
_{L^{2}\left(  K^{\ast}\right)  }+T_{k}\left(  1+\lambda\right)  \left\Vert
p-p_{k}\right\Vert _{L^{2}\left(  K\right)  }.%
\]
Hence, the claim follows from \eqref{extapprox} and $1\leq T_{k}\left(  1+\lambda\right)  $.%
\end{proof}

\printbibliography
\end{document}